\magnification1000
\tolerance 3000
\hsize = 14.5truecm
\vsize= 22truecm
\hoffset-0.2truecm
\voffset-0.5truecm
\overfullrule=0pt
\parindent=10pt

\font\bmV=cmmib10 at 5pt

\font\bmVII=cmmib10 at 7pt

\font\bmX=cmmib10
\textfont9=\bmX \scriptfont9=\bmVII \scriptscriptfont9=\bmV
\def\bm{\fam9\bmX}

\font\cmdunhX=cmr10
%\font\sb=MSBM10								
\font\sc=cmr10								
\font\scXII=cmcsc10 scaled\magstep1
\font\VII=cmr7 
\font\VIII=cmr8 
%\font\saV=msam5
%\font\saVII=msam7
%\font\saX=msam10
%\textfont10=\saX \scriptfont10=\saVII
%\scriptscriptfont10=\saV
%\def\sa{\fam10\saX}

\def\bg{\bigskip\goodbreak}
\def\bgni{\bigskip\goodbreak\noindent}

\def\mgni{\medskip\goodbreak\noindent}

\def\sgni{\smallskip\goodbreak\noindent}
\def\bnni{\bigskip\nobreak\noindent}
\let\ni=\noindent
\def\ie{{\it i.e.}}
\def\iff{if and only if }
\def\resp{\hbox{\it resp. }\ignorespaces}
\def\Block{\hbox{\vrule width 7pt height 7pt depth 0pt}}

\newcount\SecNo
\SecNo=0
\newcount\PropNo
\newcount\EqNo

% COMMENT: \Displaylines and \Eqalign are similar to TeX's
% standard \displaylines and \eqalign, excepted that each line
% can receive an equation number, by using a second & sign.
% The boolean \Displaylines is true inside a \Displaylines or
% \Eqalign group.

\newif\ifDisplaylines\Displaylinesfalse
{\catcode`\@=11
\gdef\Displaylines#1{\Displaylinestrue\displ@y\halign{
\hbox to\displaywidth{$\@lign\hfil\displaystyle##\hfil$}%
&\llap{$##$}\crcr#1\crcr}}
\gdef\Eqalign#1{\Displaylinestrue\displ@y \tabskip=\centering
\halign to \displaywidth{\hfil$\@lign\displaystyle{##}$\tabskip=0pt
&$\@lign\displaystyle{{}##}$\hfil\tabskip=\centering
&\llap{$\@lign##$}\tabskip=0pt\crcr#1\crcr}}}

% COMMENT: cross-reference system. To give give a name or a
% number to a formula or a proposition. The command \Ç54È
% causes the content of \box0 to be stored in \box54, and then
% Ç54È lets this content be printed. Similarly \Ç\nameÈ does 
% the same with a box named \name, and Ç\nameÈ lets it be
% printed. These are the only possible parameters, namely a
% number (between 1 and 255, preferably between 100 and 255),
% or a command name, i.e., \ followed by letters.

\catcode`\Ç=\active
{\catcode`\@=11
\gdef\Ç#1#2È{\ifcat\noexpand#1\noexpand\%
\global\alloc@4\box\chardef\insc@unt{#1#2}\fi%
\global\setbox#1#2\box0\ignorespaces}}
\defÇ#1È#2{\leavevmode\copy#1\ifcat#2a\ \fi#2}

\def\Sec#1\par{\vskip0pt plus.2\vsize
\penalty-250\vskip0pt plus-.2\vsize 
\advance\SecNo by 1\global\PropNo=0\global\EqNo=0
\vskip30pt\centerline{\scXII\the\SecNo. #1}}

\def\Parag#1\par{\vskip0pt plus.1\vsize
\penalty-250\vskip0pt plus -.1\vsize
\bigbreak\ni{\bf #1}}

\def\Def{\bigbreak\ni{\bf Definition. }}
\def\Not{\bigbreak\ni{\bf Notation. }}

\def\Conj{\global\advance\PropNo by 1%
\setbox0=\hbox{\the\SecNo.\the\PropNo}%
\bigbreak\ni{\bf Conjecture \the\SecNo.\the\PropNo.} \ignorespaces}

\def\Cor{\global\advance\PropNo by 1%
\setbox0=\hbox{\the\SecNo.\the\PropNo}%
\bigbreak\ni{\bf Corollary \the\SecNo.\the\PropNo.} \ignorespaces}

\def\Ex{\global\advance\PropNo by 1%
\setbox0=\hbox{\the\SecNo.\the\PropNo}%
\bigbreak\ni{\bf Example \the\SecNo.\the\PropNo.} \ignorespaces}

\def\Lem{\global\advance\PropNo by 1%
\setbox0=\hbox{\the\SecNo.\the\PropNo}%
\bigbreak\ni{\bf Lemma \the\SecNo.\the\PropNo.} \ignorespaces}

\def\Prop{\global\advance\PropNo by 1%
\setbox0=\hbox{\the\SecNo.\the\PropNo}%
\bigbreak\ni{\bf Proposition \the\SecNo.\the\PropNo.} \ignorespaces}

\def\Pr{\bigbreak\ni{\it Proof. }}
\def\EndPr\par{\hfill\nobreak$\square$\par}

\def\Eq{\global\advance\EqNo by 1
\global\setbox0=\hbox{\rm(\the\SecNo.\the\EqNo)}
\ifDisplaylines&\copy0\else\eqno{\copy0}\fi}

\def\EndPr\par{{\unskip\nobreak\hfil\penalty50%
\hskip1em\hbox{}\nobreak\Block
\parfillskip=0pt\finalhyphendemerits=0\par}}

%%%%%%%%%%%%%%%%%%%%%%%%%%%%%%%%%%%%%%%%
% 9. Basic Postscript Commands
%%%%%%%%%%%%%%%%%%%%%%%%%%%%%%%%%%%%%%%%

									% COMMENT: These TeX commands code for postcript
         					% commands and are used to describe simple drawings.
                    % All measures are in millimeters. The recommended
                    % syntax is to use a graphic box for each picture. Note
                    % that there is no PS previewer on Macintosh, and therefore
                    % the result of PS commands will not appear on the screen.
                    % However \Write is NOT a PS command.
                    % A typical example, where an external picture called "figure1"
                    % is inserted (from the file ....pictures)
                    %         $$\OpenGraphicBox width 65mm height 30mm depth 0mm;
                    %         \leavevmode\special{picture figure1 scaled 500} 
                    %         \Translate(0, 0);
                    %         \Write(35, 31, $i$);
                    %         \Write(42, 31, $i+1$);
                    %         \Write(-7, 25, $\s_i$);
                    %         \Write(-8, 10, $\s_i^{-1}$);
                    %         \CloseGraphicBox$$

\newdimen\Dwidth
\newdimen\Dheight
\newdimen\Ddepth

\def\OpenGraphicBox width#1 height#2 depth#3;
{\dimen100=\parindent\parindent=0pt
\Dwidth=#1\Dheight=#2\Ddepth=#3
\dimen104=\Dheight\advance\dimen104 by \Ddepth
\catcode`\;=14\leavevmode
\raise-\Ddepth\vbox to\dimen104\bgroup\hsize=\Dwidth
\vglue\Dheight}

\def\Translate(#1, #2){\dimen105=-#2mm\vglue\dimen105\hskip#1mm}

\def\CloseGraphicBox{\vfill\egroup\parindent=\dimen100\catcode`\;=12}

\def\Write(#1, #2, #3){\rlap{\smash{\raise#2mm\hbox{\kern#1mm #3}}}}

\def\DoASegment(#1, #2, #3, #4){
}

%%%%%%%%%%%%%%%%%%%%%%%%%%%%%%%
% Specific macros
%%%%%%%%%%%%%%%%%%%%%%%%%%%%%%%

\catcode`\¥=\active\def¥{\mathord\bullet}		% right Polish notation

\let\a=\alpha
\def\AA{{\bm A}}
\def\at(#1, #2){(#1)#2}

\let\b=\beta

		% braid group
		% braid monoid

\def\cl#1{\overline{\vrule height 5.5pt width 0pt #1}}
\def\comp{{\scriptscriptstyle\circ}}

\mathchardef\Delta="7101
\def\D#1{\Delta_{#1}}
\let\der=\partial
\def\dil{{\rm dil}}
\def\dR{\mathord{\backslash}}

\let\e=\varepsilon
\def\ea{{/\mskip -9mu o}}  % empty address 
\def\eCD{=_{\!\scriptscriptstyle C\!D}}
\let\eq=\equiv
\def\eqp{\equiv^{\scriptscriptstyle +}}		% pos equiv

\def\fCD{f_{\scriptscriptstyle \!C\!D}}  % complement for CD

\let\g=\gamma
\def\GCD{G_{\scriptscriptstyle \!C\!D}}  % group of CD
  % group of CD
\def\GGCD{{\cal G}_{\scriptscriptstyle \!C\!D}}  % geometry monoid of CD
\def\GGCDp{{\cal G}_{\scriptscriptstyle \!C\!D}^+}  % positive geometry monoid of CD
  % geometry monoid of I

\def\htR{{\rm ht}_{\!\scriptscriptstyle R}}		% right height

\def\ii{^{-1}}
\def\ince{\subseteq}

\def\Left{{\rm left}}
		% addresses

\def\MCD{M_{\scriptscriptstyle \!C\!D}}  % group of CD

\def\n{\nu}
\def\NN{{\bf N}}

\def\oa{\mathchoice{\cdot}{\mathord{\cdot}}{\mathord{\cdot}}{\mathord{\cdot}}}  % operation for addresses
\def\og{\,}  % operation in GCD
\def\om{\cdot}
  % reversed composition
\def\OP#1{\hbox{\VIII CD}^{}_{\!#1}}
\def\OPi#1{\hbox{\VIII CD}^{-1}_{\!#1}}
\def\op{*}  % operation for CD-systems

\def\opp{\mathbin{\raise3pt\hbox{$\scriptscriptstyle\wedge$}}}
\def \ot{\mathord{\cdot}}  % operations for terms

  % operation on words

\let\p=\pi
\let\pp=\dots
\def\pref{\mathrel{\raise1pt\hbox{$\scriptstyle\subset$}}}

\def\q(#1){\chi_{#1}}
\def\qs(#1){\chi_{#1}^*}

\def\RCD{R_{\scriptscriptstyle \!C\!D}}  % CD-relations

\def\sh{{\rm sh}}
\def\size{{\rm size}}

\def\Ti{T_{\!\infty}}
\def\tL#1{{t^{\scriptscriptstyle L}_{#1}}}  % left term in the identity

\def\tR#1{{t^{\scriptscriptstyle R}_{#1}}}  % right term in the identity

  % reversing relation
 % lcm

\def\wL(#1, #2){{#1}_{(#2)}}

%%%%%%%%%%%%%%%%%%%%%%%%%%%%%%%%%%%%%%%%
% List processing
%%%%%%%%%%%%%%%%%%%%%%%%%%%%%%%%%%%%%%%%

								% COMMENT: \rank(#1#2) writes the rank of the first
								% occurrence of the token #1 in the token list #2 (0 if #1 
								% does not occur.

\newif\iffound
\def\rank#1#2{{\foundfalse\count0=0 \getrank#1#2\end
\iffound\number\count0\else0\fi}} 
\def\getrank#1#2{\ifx#2\end\def\next#1{\relax}%
\else\iffound\else\advance\count0 by1\fi\let\next=\getrank
\ifx#2#1\foundtrue\fi\fi\next#1}

								% COMMENT: \List{ \a \b \c .... } causes \\\x to print the
								% rank of \x in the above list.  This strange macro is very
								% useful for references: one declares a list of references
								% at the beginning of the paper, and then typing a
								% reference with \\ in front gives the number of that
								% reference. In particular this automatically changes the
								% numbers when references are added or deleted.       					

\def\List#1{\def\\##1{\def##1{##1}\rank{##1}{#1}}}

\def\Ref#1; #2; #3; #4\par{\sgni\item{ [#1]}{\sc #2}, {\sl #3}, #4\par}
\def\Reff#1; #2; #3; #4; #5; #6; #7\par{\sgni\item{ [#1]}{\sc #2}, {\sl
#3}, #4 {\bf #5} (#6) #7\par}

\List{\Bri \CFP \Dew \Dfb \Dfg \Dfn \Dgd \Gar \Lra \Lvb \Mac \Sta}

%%%%%%%%%%%%%%%%%%%%%%%%%%%%%%%%%%%%%%%%%
%%%%%%%%%%%%%%%%%%%%%%%%%%%%%%%%%%%%%%%%%

\line{\hfill \VII 01.03}

\vskip 1 true cm 

\centerline{\cmdunhX STUDY OF AN IDENTITY}

\bg\centerline{{\sc Patrick DEHORNOY}}

\vskip 1 true cm 

{\narrower\ni{\bf Abstract}. We solve the word problem of the
identity $x(yz) = (xy)(yz)$ by investigating a certain group
describing the geometry of that identity. We also construct a
concrete realization of the free system of rank~$1$ relative to
the above identity.

\mgni Key words: word problem, free algebras, non-associative
binary operation.

\mgni MSC 2000: 03D40, 08B20, 20N02.\par}

\vskip 1 true cm

\bnni When we are given an algebraic identity~$I$ (or a family
of algebraic identities), two questions arise naturally, namely
solving the word problem of~$I$, \ie, describing an algorithm
recognizing whether two terms are forced to be equal by~$I$,
and constructing concrete realizations for the free systems in
the equational variety defined by~$I$---and, more generally, 
constructing concrete examples of systems satisfying~$I$.
Of course, answering such questions depends on the considered
identity in an essential way, and it seems hopeless to find a
uniform method that works for all identities, or, even, for a
wide class.

Due to its connection with iterations of elementary embeddings
in set theory [\\\Lvb] and with braid groups in low dimensional
topology [\\\Bri], [\\\Dgd], the left self-distributivity identity $x(yz) =
(xy)(xz)$ has received some attention in the past decade, and, in
particular, the above mentioned questions have been solved by
introducing a specific monoid that captures some    
properties of this identity [\\\Dfb], [\\\Dfn], and which turns out
to be connected with Artin's braid groups. 

A similar geometry monoid can be associated with associativity
[\\\Dfg]. In the latter case, the monoid is essentially
R.~Thompson's group~$F$ [\\\CFP], and it is closely connected
with the well known Mac~Lane--Stasheff pentagon [\\\Mac],
[\\\Sta]. Of course, solving the word problem and constructing
realizations of free systems, \ie, of free semigroups, is trivial
here.

The geometry monoid exists for every identity, and, more
generally, for every family of identities [\\\Dew]. In the most
general case, the monoid is a complicated object, of
which we have no control, and it is presumably of little help for
solving the word problem. Actually, most of the details
in~[\\\Dfb] may seem to relie on the specific properties of left
self-distributivity, making it unclear that the method can be
applied to other identities beyond the more or less trivial
case of associativity.

The aim of this paper is to show that the above mentioned
scheme does apply to other identities, yet the technical
details heavily depend on the considered identity. Here, we shall
consider
$$x(yz) = (xy)(yz), \eqno{(CD)}$$
which can be called {\it central duplication} as it consists in
duplicating the central factor~$y$. Identity~$(CD)$ has probably 
never been investigated so far, and it has probably little
interest in itself, but it should be clear that the subject of the
paper is not really that particular identity, but rather the
method we use for studying it, namely investigating the
corresponding geometry monoid.

The results we prove are:

\bgni{\bf Proposition.}
{\sl (i) The word problem of Identity~$(CD)$ is decidable, even
primitive recursive.

(ii) Let $G$ be the group $\langle \{ g_\a \, ; \, \a \in \AA\} ~;~
\RCD \rangle$, where $\AA$ is the set of all finite sequences
of~$0$'s and~$1$'s, and $\RCD$ is an effective list of relations
given in Lemma~1.3 below; let $G_0$ be the subgroup of~$G$
generated by all~$g_{0\a}$'s, and let $\sh_1$ be the
endomorphism of~$G$ that maps~$g_\a$ to~$g_{1\a}$ for
every~$\a$. Then the operation~$*$ defined on~$G$ by
$$a * b = a \om \sh_1(b) \om g_\ea \om \sh_1(b\ii)$$
induces a well defined operation on the homogeneous set~$G_0
\backslash G$; the latter operation satisfies Identity~$(CD)$,
and every monogenic subsystem of~$G_0 \backslash G$ is a free
CD-system.}

\bgni The paper is organized as follows. In Sec.$\,$1, we
introduce the geometry monoid~$\GGCD$ associated with
Identity~$(CD)$, and we establish a list of relations holding
in~$\GGCD$ called $CD$-relations. In Sec.$\,$2, we study
CD-relations from an algebraic point of view, and, in particular,
we show that the group~$\GCD$ for which CD-relations make a
presentation is a group of fractions. In Sec.$\,$3, we introduce
the blueprint of a term, which is our main tool for constructing
a binary operation satisfying a prescribed identity, here~$(CD)$.
Finally, in Sec.$\,$4, we prove the decidability of the word
problem of~$(CD)$ by using the blueprint to translate the
abstract properties of terms into concrete properties in the
group~$\GCD$.

\Sec The geometry monoid

\bnni A set equipped with a binary operation satisfying
Identity~$(CD)$ will be called a {\it CD-system}. We fix an infinite
sequence of variables $x_1$, $x_2$, \pp, and, for $1 \le n \le
\infty$, we use $T_n$ for the set of all well formed terms
constructed using $x_1$, \pp, $x_n$ and a single binary
operator. We define $\eCD$ to be the congruence relation
on~$T_n$ generated by all pairs $(t_0 \ot(t_1 \ot t_2 ), (t_0 \ot
t_1 ) \ot(t_1 \ot t_2 ))$. The quotient system $T_n / \!\! \eCD$
is the free CD-system of rank~$n$ based on~$x_1$,
\pp, $x_n$.

In order to specify geometric features precisely, we associate
with every term a finite binary tree whose leaves are labelled
with variables: if $t$ is the variable~$x$, the tree associated
with~$t$ consists of a single node labelled~$x$, while, for $t = t_0
\ot t_1 $, the tree associated with~$t$ has a root with
two immediate successors, namely a left one which is (the tree
associated with)~$t_0 $, and a right one which is (the tree
associated with)~$t_1 $. For instance, the tree associated with
$x_2 \ot (x_1 \ot x_3)$ is
$\OpenGraphicBox width 14.00mm height5.00mm
depth4.00mm;
\Translate(0, 2.00);
\Write(0.00, -2.00, $x_2$);
\Write(5.00, -4.00, $x_1$);
\Write(10.00, -4.00, $x_3$);
\DoASegment(9.50, 1.00, 7.00, -1.00);
\DoASegment(9.50, 1.00, 12.00, -1.00);
\DoASegment(5.75, 3.00, 2.00, 1.00);
\DoASegment(5.75, 3.00, 9.50, 1.00);
\CloseGraphicBox$.
We use finite sequences of~$0$'s and~$1$'s as addresses for
the nodes in such trees, starting with an empty address~$\ea$
for the root, and using $0$ and $1$ for going to the left
and to the right respectively. In this way, for each term~$t$, we
can speak of the $\a$-subterm of~$t$ for $\a$ a sufficiently short
address: for instance, the $0$-subterm (or left subterm) of~$t$
exists \iff $t$ is not a variable, and it is~$t_0$ for $t = t_0 t_1$.

\Def 
For every address~$\a$ in~$\AA$, we denote by~$\OP\a$ the
partial operator on~$\Ti$ that maps every term~$t$ with a well
defined $\a$-subterm of the form~$s_0 \ot (s_1 \ot s_2)$
to the term denoted~$\at(t, \a)$ obtained from~$t$ by replacing
the  $\a$-subterm with~$(s_0 \ot s_1) \ot (s_1 \ot s_2)$.

\bgni Thus, applying the operator~$\OP\a$ means applying
Identity~$(CD)$ in the expanding direction to the subterm with
address~$\a$. Notice that, for every~$\a$, $\OP\a$ is an
injective partial mapping on~$\Ti$, and its inverse is the
symmetric operator~$\OPi\a$ corresponding to applying~$(CD)$
in the other direction.

\Def 
The {\it geometry monoid}~$\GGCD$ of Identity~$(CD)$ is
defined to be the monoid generated by all partial
operators~$\OP\a$ and $\OPi\a$ using composition; the
submonoid of~$\GGCD$ generated by the operators~$\OP\a$
alone is denoted by~$\GGCDp$.

\bgni By construction of the congruence~$\eCD$, we have:

\Prop\Ç\CharacterizationÈ
{\sl For all terms~$t$, $t'$ in~$\Ti$, the following
are equivalent:

(i) The terms~$t$, $t'$ are CD-equivalent, \ie, $t \eCD t'$ holds;

(ii) Some element of~$\GGCD$ maps~$t$ to~$t'$.} 

\bgni By definition, every element in~$\GGCD$ is a finite
product of operators~$\OP\a$ and~$\OPi\a$. Such a product can
be specified by a word over the alphabet~$\AA \cup \AA\ii$,
where $\AA\ii$ consists of of a formal inverse~$\a\ii$ for each
address~$\a$. To this end, we define $\OP{\a\ii} = \OPi\a$, and
$\OP{uv} = \OP v \comp \OP u$---as the elements of~$\GGCD$ act
on terms on the right, it is convenient to use reversed
composition. Extending the previous notation, we write $t'
= (t)w$ when $t'$ is the image of~$t$ under~$\OP w$.
We use $\AA^*$ for the set of all words on~$\AA$, \ie, the free
monoid generated by~$\AA$, and $(\AA \cup \AA\ii)^*$ for the
set of all words on~$\AA \cup \AA\ii$. We use $\e$ for the
empty word, and define $\OP\e$ to be the identity mapping
on~$\Ti$.

The operators~$\OP w$ can be described using term unification
techniques. Let us say that a term~$t$ in~$\Ti$ is {\it canonical}
if the variables of~$t$ make an initial segment of~$(x_1, x_2,
\pp)$ when enumerated from left to right skipping repetitions;
let us say that the pair of terms~$(t_0 , t'_0)$ is an instance of
the pair~$(t, t')$ if there exists a substitution~$h$ satisfying $t_0
= h(t)$ and $t'_0 = h(t')$; finally, let us say that a term~$t$ is
{\it injective} if no variable occurs twice or more in~$t$. {\it
Mutatis mutandis}, the results of [\\\Dgd, Chap.VII] give:

\Prop\Ç\DomainÈ
{\sl (i) For every word~$w$ on~$\AA \cup \AA\ii$, either the
operator~$\OP w$ is empty, or there exists a unique pair of
CD-equivalent canonical terms $(\tL w, \tR w)$ such that $\OP
w$ maps~$t$ to~$t'$ \iff the pair $(t, t')$ is an instance of~$(\tL
w, \tR w)$. 

(ii) For every word~$u$ on~$\AA$, the operator
$\OP u$ is nonempty, and the term~$\tL u$ is injective.}

\bgni We look now for a presentation of the
monoids~$\GGCD$ and~$\GGCDp$. As in the case of left
self-distributivity [\\\Dgd] and of associativity [\\\Dfg] , we
consider relations in~$\GGCDp$ of the special type $\pp \comp
\OP\a =
\pp \comp
\OP\b$, \ie, for each pair of distinct addresses~$(\a, \b)$, we look
for possible finite sequences of addresses~$u$, $v$ satisfying
$\OP{\a \oa u} = \OP{\b \oa v}$.

\Lem\Ç\CDRelationsÈ
{\sl Let us say that a pair of words on~$\AA \cup
\AA\ii$ is a {\it CD-relation} if it is of one of the following
types:
$$\Eqalign{
\g0\a \oa \g1\b ~&,~ \g1\b \oa \g0\a 
& \hbox{(type $\perp$)}\cr
\g0\a \oa \g ~&,~ \g \oa \g00\a
& \hbox{(type $0$)}\cr
\g10\a \oa \g ~&,~ \g \oa \g01\a \oa \g10\a
& \hbox{(type $10$)}\cr
\g11\a \oa \g ~&,~ \g \oa \g11\a
& \hbox{(type $11$)}\cr
\g1 \oa \g \oa \g0 ~&,~ \g \oa \g1 \oa \g
& \hbox{(type $1$)}\cr
}$$
Then we have $\OP w = \OP{w'}$ for every CD-relation~$(w,
w')$.}

\Pr Type~$\perp$ relations are trivial. For types~$0$, $10$,
and~$11$, we observe that, when $\OP\g$ maps~$t$ to~$t'$, then
the $\g0\a$-subterm of~$t$ (if it exists) is copied to the
$\g00\a$-subterm of~$t'$,  the $\g11\a$-subterm is preserved,
and the $\g10\a$-subterm of~$t$ has two copies in~$t'$,
at~$\g01\a$ and~$\g10\a$. The last relation is less obvious,
and, in some sense, it is characteristic of the identity we
consider. Verifying the relation amounts to verifying (for $\g =
\ea$) that both $\OP{1 \oa \ea
\oa 0}$ and $\OP{\ea \oa 1 \oa \ea}$ map
$x_1 \ot (x_2 \ot (x_3 \ot x_4))$ to $((x_1 \ot x_2) \ot (x_2 \ot
x_3)) \ot ((x_2 \ot x_3) \ot (x_3 \ot x_4))$.\EndPr

\bgni At this point, we do not claim that CD-relations exhaust
all possible relations in~$\GGCD$, but we can state:

\Cor
{\sl (i) Let $\eqp$ denote the congruence on~$\AA^*$ generated
by all CD-relations. Then, for all words~$u$, $u'$ on~$\AA$, $u
\eqp u'$ implies $\OP u = \OP{u'}$.

(ii) Define $\MCD = \AA^* / \!\! \eqp$. Then $\GGCDp$ is a
quotient of~$\MCD$.}

\Sec The group~$\GCD$

\bgni An unpleasant feature with the monoid~$\GGCD$ is its
consisting of partial operators only: for every address~$\a$,
the operator~$\OP{\a \om \a\ii}$ is the identity mapping of its
domain only, and the latter is a proper subset of~$\Ti$. The
existence of words~$w$ such that $\OP w$ is empty makes it
impossible to identify all such partial identity mappings---as
is possible in the case of associativity, and, more generally, of
every identity of which both sides are injective terms. To
overcome the problem here, we introduce the group~$\GCD$ for
which CD-relations form a presentation. The leading principle in
the sequel is that $\GCD$ should resemble~$\GGCD$, and, in
particular, every notion or result about~$\GGCD$ established
using its action on terms via Identity~$(CD)$ should admit a
purely syntactic counterpart involving~$\GCD$.

\Def
We denote by $\eq$ the congruence on $(\AA \cup \AA\ii)^*$
generated by all CD-relations together with all pairs~$(\a \oa
\a\ii, \e)$ and $(\a\ii \oa \a, \e)$ with~$\a \in \AA$.  The
group~$(\AA \cup \AA\ii)^* / \! \eq$ is denoted by~$\GCD$;
for~$\a \in \AA$, the element of~$\GCD$ represented by~$\a$ is
denoted~$g_\a$. 

\bgni All subsequent results originate in the specific properties
of the group~$\GCD$, which themselves come from geometric
properties of~$(CD)$. The main technical point is that every
element of~$\GCD$ can be expressed as a right fraction of the
form~$ab\ii$, where $a$ and $b$ admit expressions where no
negative letter~$\a\ii$ occurs. This follows from the existence of
right lcm's in the monoid~$\MCD$, which will be proved now
using a uniform method called word redressing.

By definition, the CD-relations involved in the presentation of
the monoid~$\MCD$ and of the group~$\GCD$ all are of the type
$$\a \oa \pp = \b \oa \pp,$$
and, more precisely, for every pair of addresses~$(\a, \b)$,
there exists exactly one CD-relation of this type. Let us define
the mapping
$\fCD : \AA \times \AA
\rightarrow \AA^*$ by
$$\fCD(\a, \b) = \cases{
\e & for $\a = \b$, \cr
\a00\g & for $\b = \a0\g$, \cr
\a01\g \oa \a10\g & for $\b = \a10\g$, \cr
\a1 \oa \a & for $\b = \a1$, \cr
\b \oa \b0 & for $\a = \b1$, \cr
\b & in all other cases. \cr
}$$
Then $\eqp$ is the congruence on~$\AA^*$ generated by all
pairs $(\a \fCD(\a, \b), \b \fCD(\b, \a))$. A general study of those
monoids and groups with a presentation associated with a
mapping~$f$ as above can be developed along the lines of
Garside's seminal paper~[\\\Gar]. Here we extract those results
needed for our current approach. We refer to
[\\\Dgd, Chap.II] for proofs.

By construction, we have $\a f(\a , \b) \eq \b f(\b, \a)$ for
all~$\a$, $\b$, hence $\a\ii \b \eq f(\a, \b) f(\b, \a)\ii$. Let us say
that a word~$w$ on~$\AA \cup \AA\ii$ {\it redresses} to
another word~$w'$ if $w'$ is obtained from~$w$ by
repeatedly replacing some factors~$\a\ii \b$ with the
corresponding factors $f(\a, \b) f(\b, \a)\ii$. By construction,
$w$ redressing to~$w'$ implies $w \eq w'$.

The words that are terminal for word redressing are those
words of the form~$u v\ii$, where $u$ and $v$ are words
on~$\AA$. It is easy to show that every word~$w$ on~$\AA
\cup \AA\ii$ redresses to at most one word of the form~$u
v\ii$ with $u$, $v \in \AA^*$. When they exist, the latter
words will be denoted~$N(w)$ and~$D(w)$ respectively: by
definition, $w \eq N(w) D(w)\ii$ holds, and we can think
of~$N(w)$ and~$D(w)$ as the numerator and the denominator
of~$w$.

\Def
For $u$, $v$ words on~$\AA$, we define $u \dR v = N(u\ii v)$,
when the latter exists.

\bgni The operation~$\dR$ is a partial binary operation
on~$\AA^*$. By definition, we have $\a \dR \b = f(\a, \b)$ when
$\a$, $\b$ are addresses: $\dR$ is an extension of~$f$ to
arbitrary positive words.

\Lem
{\sl Assume that $u$, $v$ are words on~$\AA$ and $u \dR v$
exists. Then $v \dR u$ exists as well, and we have $u (u \dR v)
\eqp v (v \dR u)$.}

\bgni In particular, $u \dR v = v \dR u = \e$ implies $u
\eqp v$. We say that word redressing is {\it complete} when the
converse implication holds, \ie, when word redressing always
detects positive word equivalence. This need not be the case,
but we have the following effective sufficient conditions:

\Prop \Ç\CompletenessÈ  [\\\Dgd]
{\sl Let $A$ be a nonempty alphabet. Assume that $f$ is a
mapping of~$A \times A$ to~$A^*$ such that $f(x, x) = \e$
holds for every~$x$, and $f$ satisfies the following conditions:

(i) There exists a mapping~$\n : A^* \rightarrow \NN$ such that
$\n( u x f(x, y) v) = \n ( u y f(y, x) v )$ and $\n( x u) >
\n(u)$ hold for all~$x$, $y$ in~$A$, and all $u$, $v$ in~$A^*$; 

(ii) For all $x$, $y$, $z$ in~$A$, the word
$$((x \dR y) \dR (x \dR z)) \dR ((y \dR x) \dR (y \dR z))$$
exists and it is empty.

\ni Then word redressing associated with~$f$ is complete, the
monoid~$M$ associated with~$f$ admits left cancellation, and any
two elements of~$M$ that admit a common right multiple admit
a right lcm.}

\Lem
{\sl The mapping~$\fCD$ satisfies Condition~(i) of
Proposition~Ç\CompletenessÈ.}

\Pr For $t$ a term, define the size of~$t$ to be the number
of occurrences of variables in~$t$. We observe that each
operator~$\OP\a$ increases the size of the terms. We put
$$\n(u) = \size(\tR u) - \size(\tL u),$$
where $\tL u$ and $\tR u$ are the canonical terms involved in
Proposition~Ç\DomainÈ. By construction, $u \eqp u'$ implies
$\OP u = \OP{u'}$, so $\n(u)$ depends of the $\eqp$-class of~$u$
only. Assume $\a \in \AA$, and $u \in \AA^*$. By definition, we
have $\tR{\a \ot u} = 
\at(\at(\tL{\a \oa u}, \a), u)$, hence $\at(\tL{\a \oa u}, \a) =
h(\tL{u})$ and $\tR{\a \oa u} = h(\tR{u})$ for some
substitution~$h$. We deduce
$$\eqalign{\textstyle
\n(\a \oa u) 
\textstyle = \size(\tR{\a \oa u}) - \size(\tL{\a \oa u}) 
&\textstyle = \size(\tR{\a \oa u}) - \size(\at(\tL{\a \oa u}, \a)) + 
\size(\at(\tL{\a \oa u}, \a)) - \size(\tL{\a \oa u}) \cr
&\textstyle = \size(h(\tR{u})) - \size(h(\tL{u})) + 
\size(\at(\tL{\a \oa u}, \a)) - \size(\tL{\a \om u}) \cr
&> \size(h(\tR{u})) - \size(h(\tL{u})) 
\ge \size(\tR{u}) - \size(\tL{u}) = \n(u).\cr}$$
So the mapping~$\n$ satisfies the required conditions.\EndPr

\Lem
{\sl The mapping~$\fCD$ satisfies Condition~(ii) of
Proposition~Ç\CompletenessÈ.}

\Pr A priori, a lot of cases have to be considered, according to all
possible mutual positions of three addresses~$\a$, $\b$, $\g$.
However, almost all cases are automatically satisfied, as
explained in~[\\\Dfg]. Because $\ea$, $1$, and $0$ are the only
internal addresses in the two terms $x(yz)$, $(xy)(yz)$ involved
in~$(CD)$, the only non-trivial cases here are the triples
$(\ea, 1, 11)$, $(\ea, 0, 1)$, and their permuted images. In the
first case, we find
$$\displaylines{
\cases{ \ea \dR 1 = 1 \oa \ea, \cr
1 \dR \ea = \ea \oa 0, \cr} \qquad
\cases{1 \dR 11 = 11 \oa 1, \cr
11 \dR 1 = 1 \oa 10, \cr} \qquad
\cases{\ea \dR 11 = \ea, \cr
11 \dR \ea = \ea, \cr}\cr
(\ea \dR 1) \dR (\ea \dR 11) 
= (1 \oa \ea) \dR \ea = 11 \oa 1 \oa \ea, \qquad
(1 \oa \ea) \dR (1 \oa 11)
= (\ea \oa 0) \dR (11 \oa 1) = 11 \oa 1 \oa \ea, \cr
(\ea \dR 11) \dR (\ea \dR 1) 
= 11 \dR (1 \oa \ea) = 1 \oa 10 \oa \ea \oa 0, \qquad
(11 \oa \ea) \dR (11 \oa 1)
= \ea \dR (1 \oa 10) = 1 \oa \ea \oa 01 \oa 0 \oa 10,\cr
(1 \dR 11) \dR (1 \dR \ea) 
= (11 \oa 1) \dR (\ea \oa 0) = \ea \oa 0 \oa 00, \qquad
(11 \oa 1) \dR (11 \oa \ea)
= (1 \oa 10) \dR \ea = \ea \oa 0 \oa 00,\cr
%}$$
%and, finally,
%$$\displaylines{
(11 \oa 1 \oa \ea) \dR (11 \oa 1 \oa \ea) = \e, \cr
(1 \oa 10 \oa \ea \oa 0) \dR (1 \oa \ea \oa 01 \oa 0 \oa
10) = (1 \oa \ea \oa 01 \oa 0 \oa 10) \dR (1 \oa 10 \oa \ea \oa
0) = \e, \cr
(\ea \oa 0 \oa 00) \dR (\ea \oa 0 \oa 00) = \e. \cr
}$$
The verifications are similar (and simpler) in the case
of~$(\ea, 0, 1)$.\EndPr

\bgni Applying Proposition~Ç\CompletenessÈ, we deduce:

\Prop\Ç\PropCoherentÈ
{\sl Word redressing associated with~$\fCD$ is complete, the
monoid~$\MCD$ admits left cancellation, and any two elements
of~$\MCD$ that admit a common right multiple admit a right
lcm.}

\bgni We prove now that word redressing always terminates, \ie,
equivalently, that any two elements of~$\MCD$ admit a common
right multiple. The technique we use is reminiscent of Garside's
proof that any two elements in a braid monoid~$B_n^+$ admit a
common right multiple resorting to a distinguished
element~$\D n$ that is a common multiple of all generators.
Here the monoid~$\MCD$ is not of finite type, but we can
introduce some elements~$\D t$ indexed by terms which are
local counterparts to the braids~$\D n$. The intuition for
constructing the element~$\D t$ comes from the action
of~$\MCD$ on terms. 

For $t$ a term, we define the {\it right height}~$\htR(t)$
of~$t$ to be the length of the rightmost branch in~$t$ viewed as
a tree, \ie, we put $\htR(t) = 0$ for $t$ a variable, and $\htR(t_0
\ot t_1 ) = \htR(t_1 ) + 1$.

\Def
For $\a \in \AA$, we put $\a^{(p)} = \e$ for $p \le 0$, and
$\a^{(p)} = \a1^{p-1} \oa \pp \oa \a1 \oa \a$ for $p > 0$. 

\bgni If $t$ has right height~$h$, then  $\at(t,
\ea^{(p)})$ is defined exactly for~$p < h$. In this case,
assuming $t = t_0 \ot (t_1 \ot ( \pp (t_{h-1} \ot x) \pp ))$ with $x$ a
variable, we have $\at(t, \ea^{(p)} = s_0 \ot (s_1 \ot ( \pp (s_{h-1} \ot
x) \pp ))$ with $s_i = t_i \ot (t_{i+1} \ot ( \pp (t_{p-1} \ot t_p)
\pp ))$ for $ i \le p-1$ and $s_i = t_i$ for $i \ge p$.

\Not
For $w$ a word on~$\AA \cup \AA\ii$ and $\a$ an address, we
write $\a w$ for the word obtained from~$w$ by replacing each
letter~$\g^{\pm 1}$ with the letter~$(\a\g)^{\pm 1}$---not to be
confused with the word~$\a \oa w$: for~$w$ of length~$n$,
$\a w$ has length~$n$, while $\a \oa w$ has length~$n + 1$.

\Def For $t$ is a term of right height~$h$, and $\at(t,
\ea^{(h-1)}) = s_0  \ot (s_1 \ot ( \pp (s_{h-1} \ot x) \pp ))$, we
put
$$\D t = \ea^{(h-1)} \oa 0\D{s_0} \oa 10\D{s_1}
\oa \pp \oa 1^{h-1}0\D{s_{h-1}}.$$

\smallskip\ni The inductive definition of the word~$\D
t$  makes sense as, by construction, $\size(s_i) < \size(t)$
always holds. We begin with auxiliary results.

\Lem\Ç\TechIÈ
{\sl For $0 \le p \le r-2$, $0 \le q \le r-1$, and for every
word~$u$ on~$\AA$, we have 
$$\Eqalign{
1^p \oa \ea^{(r)} 
& \eqp \ea^{(r)} \oa 01^p \oa 101^{p-1} \oa \pp \oa
1^p0 \Eq\Ç129È \cr
1^q0u \oa \ea^{(r)} 
&\eqp \ea^{(r)} \oa 01^q0u \oa 101^{q-1}0u \oa \pp \oa
1^q00u, \Eq\Ç131È\cr
1^r0u \oa \ea^{(r)}
& \eqp \ea^{(r)} \oa 01^ru \oa \pp \oa 1^{r-1}01 u \oa
1^r0u. \Eq\Ç132È \cr
\ea^{(q)} \oa \ea^{(r)} 
& \eqp \ea^{(r)} \oa 0^{(q)} \oa 10^{(q-1)} \oa \pp \oa
1^{q-1}0. \Eq\Ç130È \cr
}$$

}\Pr We prove Ç129È using induction on~$p$. For $p = 0$ and $r
\ge 2$, we have
$$\ea \oa \ea^{(r)} 
= \ea \oa 11^{(r-2)} \oa 1 \oa \ea 
\eqp 11^{(r-2)} \oa \ea \oa 1 \oa \ea 
\eqp 11^{(r-2)} \oa 1 \oa \ea \oa 0 
= \ea^{(r)} \oa 0. $$
For $p > 0$, observing that $u \eqp u'$ implies $1u \eqp 1u'$, we
obtain (we mention the type of CD-relation used at each step)
$$\Eqalign{
1^p \oa \ea^{(r)}
 = 11^{p-1} \oa 1^{(r-1)} \oa \ea 
& \eqp 1^{(r-1)} \oa 101^{p-1} \oa 1101^{p-2} \oa \pp \oa
1^p0 \oa \ea \qquad\qquad &\hbox{by ind. hyp.}\cr
& \eqp 1^{(r-1)} \oa 101^{p-1} \oa \ea \oa 1101^{p-2} \oa \pp
\oa 1^p0 & (11)\cr
& \eqp 1^{(r-1)} \oa \ea \oa 01^p \oa 101^{p-1} \oa 1101^{p-2}
\oa \pp \oa 1^p0 & (10)\cr
& = \ea^{(r)} \oa \ea \oa 01^p \oa 101^{p-1} \oa
1101^{p-2} \oa \pp \oa 1^p0.\cr
}$$
We prove Ç131È using induction on~$q$. For
$q = 0$ and $r > 0$, we have
$$0u \oa \ea^{(r)} 
\eqp 1^{(r-1)} \oa 0u \oa \ea 
\eqp 1^{(r-1)} \oa \ea \oa 00u = \ea^{(r)} \oa 00u.$$
For $q > 0$, applying the induction hypothesis to
$1^{q-1}0u$ and $\ea^{(r-1)}$ and shifting all addresses
by~$1$, we find
$$\Eqalign{
1^q0u \oa \ea^{(r)} 
& = 11^{q-1}0u \oa 1^{(r-1)} \oa \ea \cr
& \eqp 1^{(r-1)} \oa 101^{q-1}0u \oa
1101^{q-2}0u \oa \pp \oa 1^q00u \oa \ea & \hbox{by ind.
hyp.}\cr 
& \eqp 1^{(r-1)} \oa 101^{q-1}0u \oa \ea\oa
1101^{q-2}0u \oa \pp \oa 1^q00u & (11) \cr
& \eqp 1^{(r-1)} \oa \ea \oa 01^q0u \oa 101^{q-1}0u \oa
1101^{q-2}0u \oa \pp \oa 1^q00u & (10) \cr
& = \ea^{(r)} \oa 01^q0u \oa 101^{q-1}0u \oa
1101^{q-2}0u \oa \pp \oa 1^q00u. \cr
}$$
We prove~Ç132È using induction on~$r \ge 0$. For
$r = 0$, Ç132È is an equality; for $r > 0$, we find
$$\Eqalign{
1^r0u \oa \ea^{(r)} 
& = 11^{r-1}0u \oa 1^{(r-1)} \oa \ea \cr
& \eqp 1^{(r-1)} \oa 101^{r-1}u \oa 1101^{r-2}u \oa
\pp \oa 1^{r-1}01u \oa 
1^r0u \oa \ea \qquad & \hbox{by ind. hyp.}\cr
& \eqp 1^{(r-1)} \oa 101^{r-1}u \oa \ea \oa 1101^{r-2}u
\oa \pp \oa 1^{r-1}01u \oa 1^r0u & (11) \cr
& \eqp 1^{(r-1)} \oa \ea \oa 01^ru \oa 101^{r-1}u \oa
\pp \oa 1^{r-1}01u \oa 1^r0u & (10) \cr
& = \ea^{(r)} \oa 01^ru \oa 101^{r-1}u \oa
\pp \oa 1^{r-1}01u \oa 1^r0u. & \cr}$$
Finally, Ç130È follows from applying Ç129È to $\ea$, $1$, \pp,
$1^{q-1}$ respectively, and gathering the factors using
$(\perp)$-relations.\EndPr

\Lem\Ç\MainÈ
{\sl Assume $t = t_0 \ot (t_1 \ot ( \pp (t_k \ot t_*) \pp ))$.
Let $t'_i = t_i \ot (t_{i+1} \ot ( \pp (t_{k-1} \ot t_k) \pp ))$ for $i
\le k$. Then there exists a word~$u$ on~$\AA$ satisfying 
$$\D t \eqp \ea^{(k)} \oa 0\D{t'_0} \oa 10\D{t'_1} \oa \pp \oa
1^k0\D{t'_k} \oa 1^{k+1}\D{t_*} \oa u. \Eq\Ç121È$$

}\Pr We use induction on the size of~$t$. Let $h = \htR(t)$. The
result is trivial if $t$ is a variable, and, more generally, for $k =
h -1$: indeed, in this case, $\D{t_*}$ is empty, and the
right hand expression in~Ç121È with $u = \e$ is the definition
of~$\D t$. Assume now $0 \le k \le h - 2$. Write $t = t_0 \ot (t_1
\ot ( \pp (t_{h-1} \ot x) \pp ))$ with $x$ a variable. Then we have
$t_* = t_{k+1} \ot ( \pp \ot (t_{h-1} \ot x) \pp )$. Let $s_i = t_i
\ot (t_{i+1} \ot ( \pp (t_{h-2} \ot t_{h-1}) \pp ))$ for $ i < h$. By
definition, we have
$$\displaylines{
\D t = \ea^{(h-1)} \oa 0\D{s_0} \oa 10\D{s_1} 
\oa \pp \oa 1^{h-1}0\D{s_{h-1}}, \cr
1^k\D{t_*} = (1^{k+1})^{(h-k-2)} \oa 1^{k+1}0\D{s_{k+1}} \oa
\pp \oa 1^{h -1}0\D{s_{h-1}}. \cr}$$ 
By construction, we have $s_i = t_i \ot (t_{i+1} \ot (\pp ( t_k \ot
s_{k+1} ) \pp ))$ and $\size(s_i) < \size(t)$ for $i \le k$, so, by
induction hypothesis, there exists a word~$u_i$ on~$\AA$
satisfying
$$\D{s_i} \eqp \ea^{(k-i)} \oa 0\D{t'_i} \oa \pp \oa
1^{k-i}0\D{t'_k} \oa 1^{k-i+1}\D{s_{k+1}} \oa u_i.$$ 
Injecting these values in~$\D t$, and using $(\perp)$-relations
to push the factors $1^i0^{(k-i)}$ to the left and the factors
$1^i0u_i$ and $1^i01^j0\D{t'_j}$ to the right,
we obtain
$$\eqalign{
\D t \eqp &\ea^{(h-1)} 
\oa 0^{(k)} 
\oa 10^{(k-1)} \oa \pp 
\oa 1^{k-1}0^{(1)} 
\oa 00\D{t'_0} 
\oa 010\D{t'_1} 
\oa 100\D{t'_1} 
\oa 0110\D{t'_2} \oa \pp
% \oa 1010\D{t'_2}
\oa 1100\D{t'_2} \cr
& \oa \pp \oa 01^k0\D{t'_k} \oa \pp
\oa 1^k00\D{t'_k} 
\oa 01^{k+1}\D{s_{k+1}} \oa \pp
\oa 1^k01\D{s_{k+1}}
\oa 1^{k+1}0\D{s_{k+1}} \cr
& \oa 1^{k+2}0\D{s_{k+2}} \oa \pp 
\oa 1^{h-1}0\D{s_{h-1}} 
\oa 0u_0 \oa 10u_1 \oa \pp \oa 1^k0u_k. \cr}$$
Applying Ç130È with $r = h-1$ and $q = k$, then Ç131È with
$r = h-1$, $q = 0$, $u = \D{t'_0}$, then $q = 1$, $u = \D{t'_1}$,
\pp, $q = k$, $u = \D{t'_k}$ successively, and, finally, Ç132È
with $r = k+1$ and $u = \D{s_{k+1}}$, we deduce
$$\Eqalign{
\D t 
& \eqp \ea^{(k)} 
\oa 0\D{t'_0} 
\oa 10\D{t'_1} \oa \pp
\oa 1^k0\D{t'_k} 
\oa (1^{k+1})^{(h-k-2)}
\oa 1^{k+1}0\D{s_{k+1}} 
\oa \ea^{(k+1)} \cr
& \hskip30mm \oa 1^{k+2}0\D{s_{k+2}} \oa \pp 
\oa 1^{h-1}0\D{s_{h-1}} 
\oa 0u_0 \oa \pp \oa 1^k0u_k\cr
& \eqp \ea^{(k)} 
\oa 0\D{t'_0} 
\oa 10\D{t'_1} \oa \pp
\oa 1^k0\D{t'_k} 
\oa (1^{k+1})^{(h-k-2)}
\oa 1^{k+1}0\D{s_{k+1}} \cr
& \hskip30mm \oa 1^{k+2}0\D{s_{k+2}} \oa \pp 
\oa 1^{h-1}0\D{s_{h-1}} 
\oa \ea^{(k+1)} 
\oa 0u_0 \oa \pp 
\oa 1^k0u_k & (11) \cr
& = \ea^{(k)} 
\oa 0\D{t'_0} 
\oa 10\D{t'_1} \oa \pp
\oa 1^k0\D{t'_k} 
\oa 1^{k+1}\D{t_*} \oa \ea^{(k+1)} 
\oa 0u_0 \oa \pp 
\oa 1^k0u_k &\Block \cr}$$

\Lem\Ç\BasicExpÈ
{\sl Assume that $(t)\a$ is defined. Then $\a \oa v
\eqp \D t$ holds for some word~$v$ on~$\AA$.}

\Pr We use induction on the length of~$\a$ as a word on~$\{0,
1\}$. Assume first $\a = \ea$. The hypothesis that $(t)\ea$ is
defined implies $\htR(t) \ge 2$. Hence $t$ can be expressed as
$t = t_0 (t_1 t_*)$. Applying Lemma~Ç\MainÈ with $k = 2$, we
obtain 
$$\D t \eqp \ea \oa 0\D{t_0 t_1 } \oa 10\D{t_1 } \oa 11\D{t_*}
\oa u$$
which begins with~$\a$ explicitly.

Assume now $\a = 0\b$. The hypothesis that $\at(t, \a)$ is
defined implies that $t$ is not a variable, so it can be written as
$t = t_0 t_*$, with $\at(t_0 , \b)$ defined. Applying
Lemma~Ç\MainÈ with $k = 1$, we obtain $u$ satisfying
$\D t \eqp 0\D{t_0 }
\oa 1\D{t_*} \oa u$, and, applying the induction hypothesis,
we obtain $v$ satisfying $\D{t_0 } \eqp \b \oa v$. We
deduce
$$\D t 
\eqp 0\D{t_0 } \oa 1\D{t_*} \oa u
\eqp 0\b \oa 0v \oa 1\D{t_*} \oa u
= \a \oa 0v \oa 1\D{t_*} \oa u,$$
which begins with~$\a$ explicitly. The argument is similar
for~$\a = 1\b$, as $0\D{t_0 }$ and $1\D{t_*}$ commute up to
$\eqp$-equivalence.\EndPr

\Lem\Ç\IncreasingÈ
{\sl Assume that $u$ is a word on~$\AA$ and $\at(t, u)$ is
defined, say $t' = \at(t, u)$. Then we have $u \oa \D{t'} \eqp
\D{t} \oa u'$ for some word~$u'$ on~$\AA$.}

\Pr We prove using induction on~$n$ that, for every word~$u$
of length~$n$, the result is true for every term~$t$ such that
$\at(t, u)$ exists. For $n = 0$, \ie, for $u = \e$, the property is
obvious. Assume $n = 1$. Then $u$ consists of a single address,
say~$\a$. We prove the result using induction on the size
of~$t$. Let $h = \htR(t)$, $t = t_0 \ot (t_1 \ot (\pp (t_{h-1} \ot x)
\pp))$, and $\at(t, \ea^{(h-1)}) = s_0 \ot (s_1 \ot (\pp (s_{h-1}
\ot x) \pp))$. Assume $t' = \at(t, \a)$. We write similarly $t' = t'_0
\ot (t'_1 \ot (\pp (t'_{h-1} \ot x) \pp))$, and $\at(t', \ea^{(h-1)})
= s'_0 \ot (s'_1 \ot (\pp (s'_{h-1} \ot x) \pp))$. We distinguish
four cases according to~$\a$.

Assume first $\a = 1^p$ with $0 \le p \le h-3$. Then we have
$t'_i = t_i$ for $i \not= p$, and $t'_p = t_p \ot t_{p+1}$. A direct
computation gives $s'_i = \at(s_i, 1^{p-i})$ for $i \le p$, and
$s'_i = s_i$ for $i > p$. For $i \le p$, the induction hypothesis
gives a word~$u'_i$ on~$\AA$ satisfying
$1^{p-i} \oa \D{s'_i} \eqp \D{s_i} \oa u'_i$. We obtain
$$\Eqalign{
\a \oa \D{t'}
&= 1^p \oa \ea^{(h-1)} \oa 0\D{s'_0} \oa 10\D{s'_1} \oa \pp
\oa 1^{h-1}0\D{s'_{h-1}} \cr
& \eqp \ea^{(h-1)} \oa 01^p \oa \pp \oa
1^p0 \oa 0\D{s'_0} \oa 10\D{s'_1} \oa \pp \oa
1^{h-1}0\D{s'_{h-1}} & \hbox{by Ç129È\relax} \cr 
& \eqp \ea^{(h-1)} \oa 01^p \oa 0\D{s'_0} \oa \pp \oa
1^p0 \oa 1^p0\D{s'_p} \oa 1^{p+1}0\D{s_{p+1}}\oa \pp \oa
1^{h-1}0\D{s_{h-1}} & (\perp) \cr 
& \textstyle
\eqp \ea^{(h-1)} \oa 0\D{s_0} \oa 0u'_0 \oa \pp \oa
1^p0\D{s_p} \oa 1^p0u'_p \oa 1^{p+1}0\D{s_{p+1}}\oa \pp
\oa 1^{h-1}0\D{s_{h-1}} \hskip20mm & \hbox{by ind. hyp.}
\cr  &  \textstyle
\eqp \ea^{(h-1)} \oa 0\D{s_0} \oa \pp \oa
1^p0\D{s_p} \oa 1^{p+1}0\D{s_{p+1}}\oa \pp
\oa 1^{h-1}0\D{s_{h-1}} \oa 0u'_0 \oa \pp \oa 1^p0u'_p 
\qquad & (\perp) \cr 
& = \D t \oa 0u'_0 \oa \pp \oa 1^p0u'_p \cr }$$
Assume now $\a = 1^{h-2}$. We still have $s'_{h-1} = s_{h-1} (= t_{h-1})$,
but, for $i < h$, we have $s'_i = t_0 \ot (t_1 \ot ( \pp \ot ((t_{h-2}
\ot t_{h-1}) \ot t_{h-1}) \pp))$, which is not a CD-expansion
of~$s_i$. Applying Lemma~Ç\MainÈ to~$t'$ with $k = h-2$ (this is
the point), we obtain a word~$u'$ on~$\AA$ satisfying
$$\D{t'} \eqp \ea^{(h-2)} \oa 0\D{s_0} \oa 10\D{s_1} \oa \pp
\oa 1^{h-2}0\D{s_{h-2}} \oa 1^{h-1}\D{t_{h-1} \ot x} \oa u'.$$
By definition, we have $\D{t_{h-1} \ot x} = 0\D{t_{h-1}}$ and $s_{h-1} =
t_{h-1}$, so we deduce
$$\a \oa \D{t'} 
\eqp 1^{h-2} \oa \ea^{(h-2)} \oa 0\D{s_0} \oa 10\D{s_1} \oa
\pp \oa 1^{h-2}0\D{s_{h-2}} \oa 1^{h-1}0\D{s_{h-1}} \oa u'
 = \D{t} \oa u'.$$
 Assume now $\a = 1^p0\b$ with $0 \le p \le h-2$. With the
same notations, we have $t'_i = t_i$ for $i \not= p$, and $t'_p =
\at(t_p, \b)$. We deduce $s'_i = \at(s_i, 1^{p-i}0\b)$ for $i \le
p$, $s'_i = s_i$ for $i > p$. For $i \le p$, the induction hypothesis
gives  a word~$u'_i$ satisfying $1^{p-i}0\b \oa \D{s'_i} \eqp
\D{s_i} \oa u'_i$, and we find
$$\Eqalign{
\a \oa \D{t'}
& = 1^p0\b \oa \ea^{(h-1)} \oa 0\D{s'_0} \oa \pp \oa
1^p0\D{s'_p} \oa \pp \oa 1^{h-1}0\D{s'_{h-1}} \cr
& \eqp \ea^{(h-1)} \oa 01^p0\b \oa \pp \oa 1^p00\b \oa
0\D{s'_0}\oa \pp \oa 1^p0\D{s'_p} \oa \pp \oa
1^{h-1}0\D{s'_{h-1}} \qquad & \hbox{by Ç131È\relax} \cr
& \eqp \ea^{(h-1)} \oa 01^{p-1}0\b \oa 0\D{s'_1}\oa \pp \oa
1^p00\b \oa1^p0\D{s'_p} \oa \pp \oa
1^{h-1}0\D{s'_{h-1}} \qquad & (\perp) \cr 
& \eqp \ea^{(h-1)} \oa 0\D{s_0} \oa 0u'_0 \oa \pp \oa
1^p0\D{s_p} \oa 1^p0u'_p \oa \pp \oa
1^{h-1}0\D{s_{h-1}} \qquad & \hbox{by ind. hyp.} \cr 
& \eqp \ea^{(h-1)} \oa 0\D{s_0} \oa \pp \oa
1^p0\D{s_p} \oa \pp \oa 1^{h-1}0\D{s_{h-1}} \oa 0u'_0 \oa \pp
\oa 1^p0u'_p & (\perp) \cr 
& = \D{t} \oa 0u'_0 \oa \pp \oa 1^p0u'_p. \cr}$$
Finally, for $\a = 1^{h-1}0\b$, we have $t'_i = t_i$ for $i
< h$ and $t'_{h-1} = \at(t_{h-1}, \b)$, hence $s'_i = \at(s_i, 1^{h-i}\b)$.
The computation is similar to the previous one, using~Ç132È
instead of~Ç131È.

Assume now $n \ge 2$. Write $u = u_1 \oa u_2$ with
$\lg(u_1), \lg(u_2) < n$. Applying the induction
hypothesis to~$u_1$ and $u_2$, we find words~$u'_1$,
$u'_2$ satisfying
$$u \oa \D{t'} 
= u_1 \oa u_2 \oa \D{\at({\at(t, u_1)}, u_2)} 
\eqp u_1 \oa \D{\at(t, u_1)} \oa u'_2 
\eqp \D{t} \oa u'_1 \oa u'_2. \eqno{\Block}$$

\Def
For $t$ a term, we put $\der t = \at(t, \D t)$---which makes
sense, as an immediate induction shows that every term~$t$
lies in the domain of the operator~$\OP{\D t}$. 

\Lem\Ç\ExpansionÈ
{\sl Assume that $u$ is a word of length~$n$ on~$\AA$, and
$\at(t, u)$ is defined. Then $u \oa v \eqp \D t \oa
\D{\der t} \oa \pp \oa \D{\der^{n-1}t}$ holds for some
word~$v$ on~$\AA$.}

\Pr We use induction on~$n$. For $n = 0$, \ie, for $u = \e$, the
result is obvious. For $n = 1$, the result is Lemma~Ç\BasicExpÈ.
Otherwise, write $u = u' \oa \a$ with $\a$ an address, and let $t'
= \at(t, u')$, which exists by hypothesis. By construction, we
have $\at(t, u) = \at(t', \a)$, hence, by Lemma~Ç\BasicExpÈ, we
have $\a \oa v' = \D{t'}$ for some~$v'$. By induction
hypothesis, there exists $u''$ satisfying
$u' \oa u'' \eqp \D t \D{\der t} \pp \D{\der^{n-2}t}$, hence
$\at(t', u'') = \der^{n-1}t$. Hence, by Lemma~Ç\IncreasingÈ, we
have $u'' \oa \D{\der^{n-1}t} = \D{t'} \oa v''$ for some~$v''$. We
find
$$u \oa v' \oa v'' \eqp u' \oa \D{t'} \oa v'' \eqp u' \oa u''
\oa \D{\der^{n-1}t} \eqp \D{t} \oa \pp \oa \D{\der^{n-2}t} \oa
\D{\der^{n-1}t}. \eqno{\Block}$$

\Prop\Ç\ConfluenceÈ
{\sl Any two elements of~$\MCD$ admit a common right
multiple.}

\Pr Assume that $u$, $v$ are words on~$\AA$. By
Proposition~Ç\DomainÈ, the terms~$\tL u$ and $\tL v$ are
injective, which implies that some substitute of~$\tL u$ is a
substitute of~$\tL v$. Hence, some term~$t$ lies both in the
domain of~$\OP u$ and~$\OP v$. Letting $n$ be the supremum of
the lengths of~$u$ and $v$, we deduce from
Lemma~Ç\ExpansionÈ the existence of two words~$u'$ , $v'$
satisfying
$$u \oa v' \eqp v \oa u' \eqp \D t \oa \D{\der t} \oa \pp \oa
\D{\der^{n-1}t}.$$
The common class of $u \oa v'$ and $v' \oa u$ in~$\MCD$ is a
common right multiple of the classes of~$u$ and~$v$
in~$\MCD$.\EndPr

\bgni Returning to word redressing in~$\AA^*$, we deduce from
the general results of~[\\\Dgd, Chap.II]:

\Prop\Ç\ComplementIIÈ
{\sl (i) Word redressing in~$(\AA \cup \AA\ii)^*$ is convergent,
\ie, for every word~$w$ on~$\AA \cup \AA\ii$, the
words~$N(w)$ and $D(w)$ exist.

(ii) For all words~$w$, $w'$ on~$\AA \cup \AA\ii$, $w \eq w'$
holds \iff we have
$$N(w) \oa v \eqp N(w') \oa v', \qquad
D(w) \oa v \eqp D(w') \oa v' \Eq\Ç140È$$
for some words~$v$, $v'$ on~$\AA$.}

\Cor
{\sl The word problem of the monoid~$\MCD$ is decidable.}

\Pr Assume that $u$, $u'$ are words on~$\AA$. By
Proposition~Ç\PropCoherentÈ, $u \eqp u'$ holds \iff
redressing~$u\ii u'$ ends with an empty word. As we know now
that redressing~$u\ii u'$ comes to an end in a finite
number of steps, this gives an effective decision method.\EndPr

\bgni Let us come back to Identity~$(CD)$. For $t$, $t'$ terms, let
us say that $t'$ is a {\it CD-expansion} of~$t$ if $t' = (t)u$ holds
for some word~$u$ on~$\AA$. If $t'$ is a CD-expansion of~$t$,
then $t'$ and $t$ are CD-equivalent, but the converse
implication is not true: going to a CD-expansion means applying
Identity~$(CD)$ in the expanding direction only. The above
results imply strong properties for the terms~$\der t$:
Lemma~Ç\IncreasingÈ implies that the operator~$\der$ is
increasing with respect to CD-expansions, \ie, that $\der t'$ is a
CD-expansion of~$\der t$ whenever $t'$ is a CD-expansion
of~$t$, and Lemma~Ç\ExpansionÈ implies that, for every term~$t$,
$\der^n t$ is a CD-expansion of every CD-expansion of~$t$
obtained by applying~$(CD)$ $n$~times at most. Finally,
Proposition~Ç\ConfluenceÈ implies the following confluence
property:

\Prop
{\sl Any two CD-equivalent terms admit a common
CD-expansion.}

\Pr As CD-equivalence is the equivalence relation generated
by the relation of being a CD-expansion, it suffices to prove that
any two CD-expansions of a term~$t$ admit a common
CD-expansion: this follows from Proposition~Ç\ConfluenceÈ
immediately.\EndPr

\bgni Actually, Proposition~Ç\ConfluenceÈ tells us a little more,
namely that, for every term~$t$, the term~$t'$ is CD-equivalent
to~$t$ \iff the term~$\der^n t$ is a CD-expansion of~$t$ for
$n$~large enough. Building on this remark, unique normal
forms with respect to CD-equivalence can be constructed along
the lines of~[\\\Dgd, Chap.VI].

\Sec The blueprint of a term

\bnni Let us address the question of constructing a monogenic
CD-system~$(S, \op)$: the question is to construct, for each
term~$t$ in~$T_1$, an interpretation of~$t$ in~$S$ in such a way
that CD-equivalent terms receive the same interpretation.  As
the only specific algebraic systems available so far are the
geometry monoid~$\GGCD$ and its abstract version~$\GCD$, we
shall start from these structures: the core of the construction
will consist in associating with every term~$t$ in~$T_1$ a
distinguished element in~$\GCD$, or, equivalently, a
distinguished word on~$\AA \cup \AA\ii$. This word arises as a
natural translation for the following property:

\Lem\Ç\AbsorptionÈ
{\sl Define $x^{[1]} = x$, and $x^{[p+1]} = x \ot x^{[p]}$ for $p \ge 1$.
Then, for every term~$t$ in~$T_1$, and for $p$ large enough, we
have
$$x^{[p+1]} \eCD t \ot x^{[p]}. \Eq\Ç100È$$

}\Pr We use induction on~$t$. For $t = x$, Ç100È is an equality.
Otherwise, assuming $t = t_0  \ot t_1 $ and using the induction
hypothesis, we obtain for $p$ large enough
$$x^{[p+1]} \eCD t_0 \ot x^{[p]} 
\eCD t_0 \ot (t_1 \ot x^{[p-1]})
\eCD (t_0 \ot t_1 ) \ot (t_1 \ot x^{[p-1]})
\eCD (t_0 \ot t_1 ) \ot x^{[p]}
= t \ot x^{[p]}. \eqno{\Block}$$

\bgni It follows from Proposition~Ç\CharacterizationÈ that, for
every term~$t$ and for every~$p$ large enough, some operator
in~$\GGCD$ maps~$x^{[p+1]}$ to~$t \ot x^{[p]}$. It suffices
to read the inductive proof of Lemma~Ç\AbsorptionÈ to obtain
an explicit description of the involved operator.

\Def
For~$t$ a term in~$T_1$,  the {\it blueprint}~$\q(t)$ of~$t$ is
the word defined by $\q(x) = \e$ and
$$\q(t) = \q(t_0 ) \om \sh_1(\q(t_1 )) \om \ea \om
\sh_1(\q(t_1 )\ii) \hbox{\qquad for $t = t_0 \ot t_1$.}
\Eq\Ç143È$$

\Prop\Ç\EffAbsÈ
{\sl For every~$t$ in~$T_1$, $\OP{\q(t)}$ maps $x^{[p+1]}$ to~$t
\ot x^{[p]}$ for $p$ large enough.}

\Pr As for Lemma~Ç\AbsorptionÈ, we use induction
on~$t$. The result is true for $t = x$. Assume $t = t_0 \ot t_1$,
and $p$ large enough. By induction hypothesis, $\OP{\q(t_0)}$
maps $x^{[p+1]}$ to~$t_0 \ot x^{[p]}$, and $\OP{\q(t_1)}$
maps $x^{[p]}$ to~$t_1 \ot x^{[p-1]}$. Hence $\OP{1\q(t_1)}$
maps $t_0 \ot x^{[p]}$ to~$t_0 \ot (t_1 \ot x^{[p-1]})$, and,
similarly, $\OP{1\q(t_1)\ii}$ maps $t \ot (t_1 \ot x^{[p-1]})$ to $t
\ot x^{[p]}$. By composing, we obtain
$$x^{[p+1]} 
~~{\buildrel\textstyle \q(t_0) \over \mapsto}~~
t_0 \ot x^{[p]} 
~~{\buildrel\textstyle 1\q(t_1) \over \mapsto}~~
t_0 \ot (t_1 \ot x^{[p-1]})
~~{\buildrel\textstyle \ea \over \mapsto}~~
(t_0 \ot t_1 ) \ot (t_1 \ot x^{[p-1]})
= t \ot (t_1 \ot x^{[p-1]}
~~{\buildrel\textstyle 1\q(t_1)\ii \over \mapsto}~~
t \ot x^{[p]}. \eqno{\Block}$$

\bgni The idea is to use the operator~$\OP{\q(t)}$, or,
rather, the image of the word~$\q(t)$ in the group~$\GCD$, as
the interpretation of the term~$t$, which leads us to introduce
the binary operation on~$\GCD$ such that the class of~$\q(t_0
\ot t_1)$ is the product of the classes of~$\q(t_0)$
and~$\q(t_1)$.

\Def
For $u$, $v$ words on~$\AA \cup \AA\ii$, we define
$$u \op v = u \om 1v \om \ea \om 1v\ii,\Eq\Ç160È$$
and we also use $\op$ for the induced binary operation
on~$\GCD$.

\bgni With this notation, $\q()$ is the homomorphism of~$T_1$
into~$((\AA \cup \AA\ii)^*, \op)$ that maps~$x$ to~$\e$. Our
plan is to start from operation~$\op$ on~$\GCD$ to construct an
operation satisfying Identity~$(CD)$. The point is that the latter
operation does not satisfy Identity~$(CD)$, but the obstruction
to its satisfying~$(CD)$ can be measured exactly. If $t$ and
$t'$ are CD-equivalent terms, their blueprints~$\q(t)$ and
$\q(t')$ need not be $\eq$-equivalent, but some operator~$\OP
w$ maps~$t$ to~$t'$, and, therefore, $\OP{0w}$ maps
$t \ot x^{[p]}$ to~$t' \ot x^{[p]}$ for every~$p$. Hence, both
$\OP{\q(t) \om 0w}$ and $\OP{\q(t')}$ map $x^{[p+1]}$ to~$t'
\ot x^{[p]}$ for $p$ large enough. If  CD-relations axiomatize the
relations in~$\GGCD$ correctly, we can therefore expect the
equivalence $\q(t) \om 0w \eq \q(t')$ to hold---which, if true,
must be verifiable by a direct computation.

\Lem\Ç\BlueprintÈ
{\sl Assume $t' = \at(t, w)$. Then we have 
$\q(t') \eq \q(t) \oa 0w$.}

\bgni For an induction on the length of~$w$ , it suffices to
prove the result when $w$ consists of a single address~$\a$.
Then, the result follows from:

\Lem
{\sl  For all words $u$, $v$, $w$ on~$\AA \cup \AA\ii$, we have
$$\Eqalign{
(u \op v) \op (v \op w) &\eq (u \op (v \op w)) \oa 0
\Eq\Ç103È \cr
(u \oa 0w) \op v 
&\eq (u \op v) \oa 00w
\Eq\Ç104È\cr
u \op (v \oa 0w) 
&\eq (u \op v) \oa 01w
\Eq\Ç105È\cr}$$

}\Pr Applying the definition of~$\op$ and CD-relations, we find
$$\Eqalign{
(u \op v) \op (v \op w)
& = u \oa 1v \oa \ea \oa 1v\ii \oa 1v \oa
11w \oa 1 \oa 11w\ii  \oa \ea \oa 11w
\oa 1\ii \oa 11w\ii \oa 1v\ii \cr
& \eq u \oa 1v \oa 11w \oa \ea \oa 1 \oa \ea 
\oa 1\ii \oa 11w\ii \oa 1v\ii & (11)\cr
& \eq u \oa 1v \oa 11w \oa 1 \oa \ea \oa 0 
\oa 1\ii \oa 11w\ii \oa 1v\ii & (1) \cr
& \eq u \oa 1v \oa 11w \oa 1 \oa \ea 
\oa 1\ii \oa 11w\ii \oa 1v\ii \oa 0 & (\perp) \cr
& \eq u \oa 1v \oa 11w \oa 1 \oa 11w\ii \oa
\ea \oa 11w \oa 1\ii \oa 11w\ii \oa 1v\ii
\oa 0 & (11) \cr 
& = (u \op (v \op w)) \oa 0. \cr
(u \oa 0w) \op v 
& = u \oa 0w \oa 1v \oa \ea \oa 1v\ii 
 \eq u \oa 1v \oa 0w \oa \ea \oa 1v\ii & (\perp) \cr
& \eq u \oa 1v \oa \ea  \oa 00w \oa 1v\ii & (0) \cr
& \eq u \oa 1v \oa \ea \oa 1v\ii \oa 00w
 = (u \op v) \oa 00w & (\perp) \cr
u \op (v \oa 0w) 
& = u \oa 1v \oa 10w \oa \ea \oa 10w\ii \oa 1v\ii \cr
& \eq u \oa 1v \oa \ea \oa 01w \oa 10w \oa 10w\ii \oa 1v\ii 
& (10) \cr
& \eq u \oa 1v \oa \ea \oa 01w \oa 1v\ii \cr
& \eq u \oa 1v \oa \ea \oa 1v\ii \oa 01w
= (u \op v) \oa 01w & \Block\cr}$$

\bgni Formula~Ç103È tells us how to obtain a binary
operation satisfying $(CD)$ from~$\op$ on~$\GCD$: it suffices
that we collapse~$g_0$. Now Ç104È and~Ç105È show that, in
order to obtain a well defined induced operation, we have to
collapse every generator~$g_{0\a}$ as well. So we have:

\Prop\Ç\OperationÈ
{\sl For every address~$\g$, let $\sh_\g$ denote the
endomorphism of~$\GCD$ induced by the address shift $\a
\mapsto \g\a$. Then the operation~$\op$ defined on~$\GCD$ by
$$a \op b = a \oa \sh_1(b) \oa g_\ea \oa \sh_1(b\ii)$$
induces a well defined operation on the coset set $\sh_0(\GCD)
\backslash \GCD$, and the latter operation satisfies
Identity~$(CD)$.}

\bgni We shall say more about the previous operation (and, in
particular, prove that it is not trivial) in the next section. We
conclude the current section with a complete description of the
connection between the group~$\GCD$ and the geometry
monoid~$\GGCD$.

Assume that $w$ and $w'$ are words on~$\AA \cup \AA\ii$ and
both $\OP w$ and $\OP{w'}$ map the term~$t$ to the term~$t'$.
Then, by Lemma~Ç\BlueprintÈ, we have
$$\q(t) \oa 0w \eq \q(t') \eq \q(t) \oa 0w',$$
which implies $0w \eq 0w'$. We observe that, if the address~$\g$
is a prefix of all addresses involved in the left term of a
CD-relation, then the same holds for the right term, and {\it
vice versa}. It follows that $\g u \eqp \g u'$ implies $u \eqp u'$
for all words~$u$, $u'$ on~$\AA$, as all intermediate words in a
sequence of elementary transformations from~$\g u$ to~$\g u'$
witnessing $\g u \eqp \g u'$ must be of the form~$\g v$. 
Now, arbitrary factors $\a \oa \a\ii$ may appear in
$\eq$-equivalences, and the same argument does not apply
to~$\eq$. It is actually true that $0w \eq 0w'$ implies $w \eq w'$,
but the proof requires a number of auxiliary results. We can
avoid the problem by resorting to an alternative blueprint.

\Def
For $t$ in~$T_1$, we define $\qs(t) = \e$ for $t = x$,
and $\qs(t) = \q(t_0 ) \oa 1\qs(t_1 )$ for $t = t_0 \ot
t_1$.

\Lem\Ç\HalfBlueprintÈ
{\sl Assume $t' = \at(t, w)$. Then we have $\qs(t') \eq \qs(t)
\oa w$.}

\Pr For an induction, it suffices to prove the result when $w$
consists of a single positive address, say~$\a$. We use induction
on the length of~$\a$ as a word on~$\{0, 1\}$. Assume first $\a =
\ea$. As $\at(t, \ea)$ exists, we can write $t = t_0 \ot
(t_1 \ot t_*)$, and we have then $t' = (t_0 \ot t_1) \ot (t_1 \ot
t_*)$. We find
$$\Eqalign{
\qs(t')
= \q(t_0 \ot t_1) \oa 1\q(t_1) \oa 11\qs(t_*) 
&= \q(t_0) \oa 1\q(t_1) \oa \ea \oa 1\q(t_1)\ii \oa 1\q(t_1)
\oa 11\qs(t_*) \cr
&\eq \q(t_0) \oa 1\q(t_1) \oa \ea \oa 11\qs(t_*) 
\eq \q(t_0) \oa 1\q(t_1) \oa 11\qs(t_*) \oa \ea 
= \qs(t) \oa \ea. \cr }$$
Assume now $\a = 0\b$. Write $t = t_0 \ot t_1$. Then we have $t'
= t'_0 \ot t_1$ with $t'_0 = \at(t_0, \b)$. Applying
Lemma~Ç\BlueprintÈ, we find
$$\qs(t')
= \q(t'_0) \oa 1\qs(t_1)
\eq \q(t_0) \oa 0\b \oa 1\qs(t_1)
\eq \q(t_0) \oa 1\qs(t_1) \oa 0\b
= \qs(t) \oa \a.$$
Assume finally $\a = 1\b$. We write $t = t_0 \ot t_1$ again. Then
we have $t' = t_0 \ot t'_1$ with $t'_1 = \at(t_1, \b)$. Applying the
induction hypothesis, we find 
$$\qs(t')
= \q(t_0) \oa 1\qs(t'_1)
\eq \q(t_0) \oa 1\qs(t_1) \oa 1\b 
= \qs(t) \oa \a.
\eqno{\Block}$$

\Lem\Ç\ReversingÈ
{\sl Assume that $w$ is a word on~$\AA \cup \AA\ii$, $w$
redresses to~$w'$, and $\at(t, w)$ is defined. Then $\at(t, w')$
is defined as well.}

\Pr It suffices to consider the case where exactly one factor
$\a\ii \oa \b$ is replaced with the corresponding factor
$\fCD(\a, \b) \oa \fCD(\b, \a)\ii$. Then we consider each possible
CD-relation. The details are easy.\EndPr

\Prop
{\sl Assume that $w$, $w'$ are words on~$\AA \cup \AA\ii$, and
the domains of~$\OP w$ and~$\OP{w'}$ are not disjoint. Then the
following are equivalent:

- ~ $\at(t, w) = \at(t, w')$ holds for at least one term~$t$;

- ~ $\at(t, w) = \at(t, w')$ holds for every term~$t$ such that 
$\at(t, w)$ and $\at(t, w')$ exist;

- ~ $w \eq w'$ holds.

\ni If $w$ and $w'$ are words on~$\AA$, $\OP w = \OP{w'}$ is
equivalent to $w \eq w'$, so $\GGCD^+$ is isomorphic to the
submonoid~$\GCD^+$ of~$\GCD$ generated by the
elements~$g_\a$.}

\Pr Assume that both $\OP w$ and $\OP{w'}$ map $t$
to~$t'$. By Lemma~Ç\HalfBlueprintÈ, we have
$$\qs(t) \om w \eq \qs(t') \eq \qs(t) \om w',$$ 
hence $w \eq {\qs(t)}\ii \oa \qs(t') \eq w'$.

Conversely, assume that $w \eq w'$ holds, and both $\at(t, w)$
and $\at(t, w')$ exist. By Lemma~Ç\ReversingÈ,  $\at(t, {N(w)
D(w)\ii})$ and
$\at(t, {N(w') D(w')\ii})$ exist. By Proposition~Ç\ComplementIIÈ,
there exists two words~$v$, $v'$ on~$\AA$ satisfying $N(w) \, v
\eqp N(w') \, v'$ and $D(w) \, v \eqp D(w') \, v'$, and we find
$$\Eqalign{
\at(t, w) 
= \at(t, {N(w) \, D(w)\ii})
& = \at(t, {N(w) \, v \, v\ii \, D(w)\ii})\cr
&= \at(t, {N(w') \, v' \, {v'}\ii \, D(w')\ii})
= \at(t, {N(w') \, D(w')\ii})
= \at(t, w').\cr}$$
If, in addition, $w$ and $w'$ are words on~$\AA$, then the
terms~$\tL{w}$ and $\tL{w'}$ are injective, and the basic
properties of term unification imply that the domains of~$\OP
w$ and $\OP{w'}$ are never disjoint: the previous results
apply, so $w \eq w'$ implies that $\OP w$ and $\OP{w'}$
agree on every term on which they are both defined. To
conclude that $\OP w$ and $\OP{w'}$ coincide, we resort
to the results of [\\\Dgd, Chapter~VII], which apply {\it
mutatis mutandis}.\EndPr

\Cor
{\sl The word problem of the group~$\GCD$ is decidable.}

\Pr Assume that $w$ is a word on~$\AA \cup \AA\ii$. Then $w
\eq \e$ is equivalent to~$N(w) \eq D(w)$, hence, by the previous
result, to~$\at(t, {N(w)}) = \at(t, {D(w)})$ for some/any
term~$t$ in the intersection of the domains of~$\OP{N(w)}$
and~$\OP{D(w)}$. Such a term~$t$ can be computed effectively
from~$w$ and~$w'$ using unification.\EndPr

\Sec Iterated left subterms

\bnni Let us consider the word problem of Identity~$(CD)$, \ie,
the problem of recognizing CD-equivalent terms. In the case of
one variable terms, Lemma~Ç\BlueprintÈ tells us that $t \eCD t'$
implies that the class of $\q(t)\ii \oa
\q(t)$ in the group~$\GCD$ belongs to the
subgroup~$\sh_0(\GCD)$. At this point, we do not know that the
previous implication is an equivalence, and we have no
effective criterion for recognizing elements of~$\sh_0(\GCD)$.
The last ingredient needed in our construction is a
preordering on~$\GCD$ enabling us to prove that a given
element of~$\GCD$ does not belong to~$\sh_0(\GCD)$. Once
again, the considered property of~$\GCD$ is the translation of
some geometric feature involving Identity~$(CD)$, namely the
action on iterated left subterms.

If $t'$ is a CD-expansion of~$t$, then some iterated left
subterm of~$t'$ is a CD-expansion of the left subterm of~$t$, as a
trivial induction shows. For $t$ a term that is not a variable, let
us denote by~$\Left(t)$ the left subterm of~$t$. The precise
statement is as follows:

\Lem\Ç\LeftExpÈ
{\sl Define~$\dil: \NN \times \AA^* \rightarrow \NN$ 
inductively by
$$\dil(i , \e) = i, 
\quad \dil(i, \a) = \cases{i+1 & for $\a = 1^p$ with $p < i$, \cr 
i & otherwise, \cr}, 
\quad \dil(i, u \oa v) = \dil(\dil(i, u), v).$$
Assume that $u$ is a word on~$\AA$, and $t' = \at(t, u)$ holds.
Then, for every~$i$ such that $\Left^i(t)$ exists,
$\Left^{\dil(i, u)}(t')$ is a CD-expansion of~$\Left^i(t)$.}

\bgni The proof is an easy induction. If $u$ and $u'$ are
$\eqp$-equivalent words on~$\AA$, the operators~$\OP u$ and
$\OP{u'}$ coincide, and we can therefore expect the
mappings~$\dil(\cdot, u)$ and~$\dil(\cdot, u')$ to coincide as
well. Once again, if true, this property must be easily verifiable.

\Lem\Ç\DilÈ
{\sl Assume $u , u' \in \AA^*$ and $u \eqp u'$. Then we have
$\dil(i, u) = \dil(i, u')$ for every~$i$.}

\Pr Consider all basic CD-relations successively.\EndPr

\bgni When we consider an word~$w$ on~$\AA \cup \AA\ii$,
the integers~$\dil(i, w)$ are no longer defined, but we can
consider the values associated with the numerator and the
denominator of~$w$. These values depend on~$w$, but their
relative position depends on the $\eq$-class of~$w$ only:

\Lem\Ç\DilCompÈ
{\sl Assume $w , w' \in (\AA \cup \AA\ii)^*$ and $w \eq w'$.
Then $\dil(1, D(w)) = \dil(1, N(w))$ (\resp
$<$, $>$) is equivalent to $\dil(1, D(w')) = \dil(1, N(w'))$ (\resp
$<$, $>$).}

\Pr By Proposition~Ç\ComplementIIÈ, there exist words $v$,
$v'$ on~$\AA$ satisfying
$N(w) \, v \eqp N(w') \, v'$ and $D(w) \, v
\eqp D(w') \, v'$. Applying the definition of~$\dil$ and
Lemma~Ç\DilÈ, we find
$$\displaylines{
\dil(\dil(1, D(w)), v) = \dil(1, D(w) v) 
= \dil(1, D(w') v') = \dil(\dil(1, D(w')), v'), \cr
\dil(\dil(1, N(w)), v) = \dil(1, N(w) v) 
= \dil(1, N(w') v') = \dil(\dil(1, N(w')), v').\cr}$$
By construction, the mappings $\dil(\cdot, v)$ and
$\dil(\cdot, v')$ are increasing, hence
$\dil(1, D(w)) = \dil(1, N(w))$ is equivalent to
$\dil(1, D(w')) = \dil(1, N(w'))$, and the same for~$<$
and~$>$. \EndPr

\Prop\Ç\OneVariableÈ
{\sl Assume that $t$, $t'$ are terms in~$T_1$. Then the following
are equivalent:

(i) The terms~$t$ and $t'$ are CD-equivalent;

(ii) We have $\dil(1, D(\q(t)\ii \oa \q(t'))) = \dil(1, N(\q(t)\ii
\oa \q(t')))$.}

\Pr Let $w = \q(t)\ii \oa \q(t')$. Assume~(i). By
Lemma~Ç\BlueprintÈ, we have $w \eq 0w_0$ for some
word~$w_0$. By construction, we have
$D(0w_0) = 0D(w_0)$ and $N(0w_0) = 0N(w_0)$, and $\dil(1, 0u)
= 1$ for every word~$u$ on~$\AA$. Hence we have $\dil(1,
D(0w_0)) = \dil(1, N(0w_0)) = 1$, which, by Lemma~Ç\DilCompÈ,
implies~$\dil(1, D(w)) = \dil(1, N(w))$.

Assume now~(ii). By Proposition~Ç\EffAbsÈ and
Lemma~Ç\ReversingÈ, we have 
$$\at(t \ot x^{[p]}, w) 
= \at(t \ot x^{[p]}, {N(w) \oa D(w)\ii})
= t' \ot x^{[p]}$$
for $p$ large enough. Let $t_0 = \at(t \ot x^{[p]}, {N(w)})$. By
construction, we have $t_0 = \at(t' \ot x^{[p]}, {D(w)})$. Let $k$
be the common value of~$\dil(1, D(w))$ and $\dil(1, N(w))$. By
lemma~Ç\LeftExpÈ, $\Left^k(t_0)$ is a CD-expansion both
of~$\Left(t \ot x^{[p]})$, \ie, of~$t$, and of
$\Left(t' \ot x^{[p]})$, \ie, of~$t'$. It follows that $t$ and $t'$ are
CD-equivalent, since they admit a common
CD-expansion.\EndPr

\Cor
{\sl The word problem of Identity~$(CD)$ restricted to one
variable terms is decidable.}

\Pr The integers~$\dil(1, N(\q(t)\ii \oa \q(t')))$
and $\dil(1, D(\q(t)\ii \oa \q(t')))$ are effectively computable.
\EndPr

\bgni Extending the solution of the word problem to the general
case turns out to be easy.

\Lem\Ç\AcyclÈ
{\sl (i) A term is never CD-equivalent to one of its proper
iterated left subterms.

(ii) Distinct terms with the same skeleton are never
CD-equivalent.}

\Pr (i) For $w$ a word on~$\AA \cup \AA\ii$, let us say that $w$
belongs to~$P_+$ (\resp $P_0$) if $\dil(1, D(w)) <
\dil(1, N(w))$ holds (\resp $=$). By Lemma~Ç\DilCompÈ, the
sets~$P_+$ and $P_0$ is saturated under~$\eq$. Assume $w_1$,
$w_2 \in P_+$. We find
$$\eqalign{
\dil(1, D(w_1w_2))
& = \dil(1, D(w_2) \, (N(w_2) \dR D(w_1))) \cr
& = \dil(\dil(1, D(w_2)) , N(w_2) \dR D(w_1)) \cr
& < \dil(\dil(1, N(w_2)) , N(w_2) \dR D(w_1)) \cr
& = \dil(1, N(w_2) \, (N(w_2) \dR D(w_1))) \cr
& = \dil(1, D(w_1) \, (D(w_1) \dR N(w_2))) \cr
& = \dil(\dil(1, D(w_1)), D(w_1) \dR N(w_2)) \cr
& < \dil(\dil(1, N(w_1)), D(w_1) \dR N(w_2)) \cr
& = \dil(1, N(w_1) \, (D(w_1) \dR N(w_2))) 
=\dil(1, N(w_1w_2)), \cr}$$
so we have $P_+ \oa P_+ \ince P_+$, and, by a similar
argument, $P_0 \oa P_+ \ince P_+$, and $P_+ \oa P_0 \ince P_+$.

Assume now that $t$ is a proper iterated left subterm of~$t'$:
this means that we have $t' = (( \pp (t \ot t_1) \ot t_2) \ot \pp)
\ot t_k$ for some terms~$t_1$,\pp, $t_k$, which, by definition,
gives a decomposition of the form
$$\q(t') = \q(t) \oa 1w_0 \oa \ea \oa 1w_1 \oa \pp \oa
1w_{k-1} \oa \ea \oa 1w_k.$$
For each~$i$, the word~$1w_i$ belongs to~$P_0$, while $\ea$
belongs to~$P_+$, since we have $\dil(1, D(\ea)) = \dil(1, \e) = 1$
and $\dil(1, N(\ea)) = \dil(1, \ea) = 2$. By the above
computations, we deduce $\q(t)\ii \oa \q(t') \in P_+$, while,
by Proposition~Ç\OneVariableÈ, $t \eCD t'$  is equivalent to
$\q(t)\ii \oa \q(t') \in P_0$.

(ii) Assume that $t$, $t'$ are distinct terms with the same
skeleton. Assume that some variable~$x$ occurs at~$\a$ in~$t$,
while $x'$ occurs at~$\a$ in~$t'$. We assume
$(x, x')$ to be the leftmost variable clash between~$t$ and~$t'$.
First, by replacing~$t$ and ~$t'$ by some CD-expansion, we can
assume that $\a$ has the form~$0^i1^j$, \ie, the clash involves
the rightmost variable in the $p$-th iterated left subterm of~$t$
and~$t'$. Let $t''$ be a common CD-expansion for~$t$ and~$t'$.
By Lemma~Ç\LeftExpÈ, we have $\Left^k(t'') \eCD \Left^i(t)$
and $\Left^{k'}(t'') \eCD \Left^i(t')$ for some~$k$, $k'$. As the
rightmost variables in~$\Left^i(t)$ and $\Left^i(t')$ are distinct,
and the righmost variable is preserved under CD-equivalence,
we deduce $k \not= k'$. Assume for instance $k > k'$. Then
$\Left^k(t'')$ is a proper iterated subterm of~$\Left^{k'}(t'')$,
hence of~$t_0$, where $t_0$ is the term obtained
from~$\Left^{k'}(t'')$ by replacing the final variable~$x'$ by~$x$.
Now $t_0$ is CD-equivalent to the term obtained
from~$\Left^i(t')$ by replacing the final variable with~$x$, and
the latter term is~$\Left^i(t)$. So $t_0$ is CD-equivalent to its
proper iterated subterm~$\Left^{k-k'}(t_0)$, contradicting~(i).
\EndPr

\Prop
{\sl The word problem of Identity~$(CD)$ is decidable, with a
primitive recursive complexity.}

\Pr Assume that $t$, $t'$ are terms in~$\Ti$. Let $t_1$ and $t'_1$
respectively be the terms in~$T_1$ obtained by replacing every
variable in~$t$ and~$t'$ with~$x_1$. We can decide $t \eCD
t'$ as follows. First, we test $t_1 \eCD t'_1$ using
Proposition~Ç\OneVariableÈ. If $t_1 \eCD t'_1$ fails, so does $t
\eCD t'$. Otherwise, $t_1$ and $t'_1$ admit a common
CD-expansion, namely $\at(t_1, u) = \at(t'_1, u')$, with
$u = N(\q(t_1)\ii \oa \q(t'_1))$ and
$u' = D(\q(t_1)\ii \oa \q(t'_1))$. Then we compare $\at(t, u)$ and
$\at(t' , u')$: these terms exist, as, for $u$ in~$\AA^*$, $\at(t, u)$
being defined only depends on the skeleton of~$t$, and they
have the same skeleton, namely the common skeleton of
$\at(t_1, u)$ and $\at(t'_1, u')$. Then $\at(t, u') = \at(t' , u')$
implies $t \eCD t'$, while, by Lemma~Ç\AcyclÈ(ii), $\at(t, u) \not=
\at(t' , u')$ implies $\at(t, u) \not\eCD \at(t' , u')$, hence $t
\not\eCD t'$.

As for complexity, we observe that, if $t$ and $t'$ have
size~$n$ at most, then the whole computation can be made
using space resources not larger than the size of the
term~$\der^{2^n} x^{[n]}$, and the latter is bounded above by a
tower of exponentials of height~$2^n$.\EndPr

\bgni If $S$ is an arbitrary binary system, we say that $a$ is a
left divisor of~$b$ if $b = ax$ holds for some~$x$, and that $a$ is
an iterated left divisor of~$b$, denoted $a \pref b$, if we have $b
= ( \pp (a x_1) \pp ) x_k$ for some~$x_1$, \pp, $x_k$ (the two
notions coincide in the case of an associative operation only).

\Prop\Ç\OrderÈ
{\sl Assume that $S$ is a free CD-system. Then iterated left
division is a partial order on~$S$. Moreover, if $S$ has rank~$1$,
this order is a linear order.}

\Pr As $\pref$ is transitive by definition, the point is to prove
that $\pref$ is irreflexive, which follows from
Lemma~Ç\AcyclÈ(i): indeed, assume that $a$ is the class of the
term~$t$; then $a \pref a$ is equivalent to the existence
of a term~$t'$ such that $t'$ is  CD-equivalent to~$t$ and $t$ is
CD-equivalent to a proper iterated left subterm of~$t'$.

Assume now that $S$ is a free CD-system of rank~$1$, and
$a$, $a' \in S$ holds. Let $t$, $t'$ be one variable terms
representing~$a$ and $a'$ respectively. Let $w = \q(t)\ii \oa
\q(t')$. Let $t_0 = \at(t \ot x^{[p]}, {N(w)})$,
$k = \dil(1, N(w))$, and $k' = \dil(1, D(w))$. As in the proof of
Proposition~Ç\OneVariableÈ, we see that $\Left^k(t_0)$ is a
CD-expansion of~$t$, while $\Left^{k'}(t_0)$ is a CD-expansion
of~$t'$, and, therefore, $\Left^k(t_0)$ represents~$a$ and
$\Left^{k'}(t_0)$ represents~$a'$. For $k = k'$, we obtain $a = a'$.
For $k > k'$, the term~$\Left^k(t_0)$ is a proper iterated left
subterm of~$\Left^{k'}(t_0)$, and we deduce $a \pref a'$.
Similarly $k < k'$ implies $a' \pref a$.\EndPr

\bgni An application of the previous results is the following
criterion for recognizing free CD-systems, which is directly
reminiscent of Laver's criterion for free LD-systems [\\\Lvb]:

\Prop\Ç\CriterionÈ
{\sl A monogenic CD-system~$S$ is free \iff left division has no
cycle in~$S$.}

\Pr By Proposition~Ç\OrderÈ, the condition is necessary.
Conversely, assume $S$ to be generated by~$g$. Let $F$ be a
free CD-system based on~$\{x\}$, and let $\p$ be the canonical
homomorphism of~$F$ onto~$S$ that maps~$x$ to~$g$. Let $a$,
$a'$ distinct elements of~$F$. By Proposition~Ç\OrderÈ, either $a
\pref a'$ or $a' \pref a$ holds. As $\pref$ is definable from the
binary operation, $\p$ preserves~$\pref$, so $\p(a)
\pref \p(a')$ or $\p(a') \pref \p(a)$ holds in~$S$. If left division
in~$S$ has no cycle, both imply $\p(a) \not= \p(a')$, $\p$ is
injective, and $S$ is isomorphic to~$F$, hence free.\EndPr

\bgni Let us come back to the CD-system~$(\sh_0(\GCD)
\backslash \GCD, \op)$ of Proposition~Ç\OperationÈ. For
simplicity, we write $G$ for $\GCD$ and $G_0$ for $\sh_0(\GCD)$
in the sequel. The operation~$\op$ on~$G_0 \backslash G$ is
defined by
$$aG_0 \op bG_0 = a \, \sh_1(b) \, g_\ea , \sh_1(b\ii)\, G_0.$$
The remaining question is whether the latter binary operation
is trivial or not: when collapsing all generators~$g_{0\a}$
in~$G$, we could have collapsed everything and obtained a
trivial quotient. Actually, we have not:

\Prop
{\sl Every monogenic subsystem of $(G_0 \backslash G,
\op)$ is free.}

\Pr Assume that $a_0G_0$, \pp, $a_kG_0$ are cosets in~$G_0 \backslash
G$ and each factor divides the next one, \ie, we have $a_iG_0
\op x_iG_0 = a_{i+1}G_0$ for some~$x_i$. This means that, for
every~$i$,  we have $(a_i \op x_i ) \oa \sh_0(y_i) = a_{i+1}$
in~$G$ for some~$y_i$. By using the definition of~$\op$ and
gathering the equalities, we obtain in$~G$ an equality of the
form
$$a_k = a_0 \, \sh_1(c_0)  \, g_\ea  \, \sh_1(c_1) \sh_0(c'_1)
 \, g_\ea \,  \pp  \, g_\ea \,  \sh_1(c_k) \sh_0(c'_k). \Eq\Ç180È$$
For $k \ge 1$, Ç180È shows that $a_1\ii a_k$ can be represented
by a word containing $k$~letters~$\ea$, and no letter~$\ea\ii$,
hence a word in the set~$P_+$ introduced in the proof of
Lemma~Ç\AcyclÈ, and, therefore, not in~$P_0$, as would
be the case if we had $a_0\ii a_k \in G_0$. So we deduce $a_kG_0
\not= a_0G_0$, \ie, $(a_0G_0, \pp, a_kG_0)$ is not a cycle for left
division in~$(G_0 \backslash G, \op)$. Proposition~Ç\CriterionÈ
then implies that every mongenic subsystem of $(G_0
\backslash G, \op)$ is free.\EndPr

\bgni{\bf Remarks.}
(i) If, for $a, b$ in~$G$, we say that $a \prec b$ (\resp $a
\simeq b$) holds if $a\ii b$ admits an expression in~$P_+$ (\resp
in~$P_0$), then ~$\prec$ is a preorder on~$G$,
and $\simeq$ is the associated equivalence relation; both are
invariant under left multiplication. The previous proof means
that $a \prec a \op b$ holds for all $a$, $b$ in~$G$, and the
preorder~$\prec$ is connected with the iterated left divisibility
order~$\pref$ on free CD-systems of rank~$1$ as follows: for
$t$, $t'$ in~$T_1$, $\cl{t} \pref \cl{t'}$ holds in~$T_1 / \!\! \eCD$
\iff $\cl{\q(t)} \prec \cl{\q(t')}$ holds in~$G$, where $\cl{t}$
denotes the $\eCD$-class of~$t$, and $\cl{w}$ the $\eq$-class
of~$w$.

(ii) If, instead of considering the cosets associated with the
subgroup~$G_0$, we consider the normal
subgroup~$\widehat{G_0}$ of~$G$ generated by~$G_0$, we
still obtain an operation satisfying~$(CD)$ on the
quotient-group~$G / \! \widehat{G_0}$---but the latter
quotient is trivial: for  every address~$\g$, the
CD-relation~$g_\g \og g_{\g1} \og g_\g = g_{\g1} \og g_\g \og
g_{\g0}$ in~$G$ implies
$g_\g \og g_{\g1}
\og g_\g = g_{\g1} \og g_\g$, hence $g_\g = 1$, in~$G / \!
\widehat{G_0}$, and  $\widehat{G_0}$ is all of~$G$. This
distinguishes $(CD)$ from left self-distributivity~$(LD)$: in the
latter case, we have a similar situation where a binary operation
satisfying~$(LD)$ exists both on a coset set~$G'_0 \backslash
G'$---where $G'$ is a certain group connected with the
geometry monoid of~$(LD)$---and on the quotient group~$G' /
\! \widehat{G'_0}$, where $\widehat{G'_0}$ is the normal
subgroup of~$G'$ generated by~$G'_0$; now $G' / \!
\widehat{G'_0}$ turns out to be Artin's braid group~$B_\infty$,
and one can deduce a simplified solution for the word problem
of~$(LD)$ by using this group and its representation in the
automorphisms of a free group [\\\Lra].  In the case of~$(CD)$,
such an indirect approach is not possible.

\bgni The study of Identity~$(CD)$ can be continued along the
lines developed for left self-distributivity in~[\\\Dgd]. As natural
examples are missing, going into details seems unnecessary. Let
us only mention that the group~$\GCD$ is an orderable group,
\ie, it can be equipped with a linear ordering compatible with
multiplication on both sides, and that one can construct
realizations for the free CD-systems of any rank by extending
the blueprints so as to generate arbitrary terms.

As it stands, the current analysis, which is reminiscent of
Henkin's proof of G\"odel's completeness theorem, relies on
three ingredients, namely the completeness of the
involved word redressing, its convergence, and the existence of a
convenient blueprint. We conjecture that the first
condition holds whenever the left term of the
considered identity is injective, \ie, no variable is
repeated. For the other conditions, no general principle
arises so far, but, in any case, the current scheme is not the
only possible one for using the geometry monoid, and we hope
for new applications of the latter in the future.

\bg\bgni\centerline{\scXII References}\bnni
\parindent=20pt

\Reff \\\Bri; E. Brieskorn; Automorphic sets and braids and
singularities; Braids, Contemporary Maths AMS; 78; 1988;
45--117.

\Reff \\\CFP; J.W. Cannon, W.J. Floyd, \& W.R. Parry;
Introductory notes on Richard Thompsons's groups; Ens.
Math.; 42; 1996; 215--257.

\Reff \\\Dew; P. Dehornoy; Structural monoids associated to
equational varieties; Proc. Amer. Math. Soc.; 117-2; 1993;
293--304.

\Reff \\\Dfb; ---; Braid groups and left
distributive operations; Trans. Amer. Math. Soc.; 345-1;
1994; 115--151. 

\Reff \\\Dfg; ---; The structure group for the
associativity identity; J. Pure Appl. Algebra; 111; 1996;
59--82.

\Reff \\\Dfn; ---; Construction of left
distributive operations and charged braids; Int. J. for
Algebra \& Computation; 10-1; 2000; 173--190.

\Ref \\\Dgd; ---; Braids and Self-Distributivity;
Progress in Math. vol. 192, Birkh\"auser, (2000).

\Reff \\\Gar; F. A. Garside; The braid group and
other groups; Quart. J. Math. Oxford; 20 {\rm No.78}; 1969;
235--254.

\Reff \\\Lra; D.M. Larue; On braid words and irreflexivity;
Algebra Univ.; 31; 1994; 104--112.

\Reff \\\Lvb; R. Laver; The left distributive law and the
freeness of an algebra of elementary embeddings;
Advances in Math.; 91-2; 1992; 209--231. 

\Reff \\\Mac; S.~Mac Lane; Natural associativity and
commutativity; Rice Univ. Studies; 49; 1963; 28--46

\Reff \\\Sta; J. Stasheff; Homotopy associativity of $H$-spaces;
Trans. Amer. Math. Soc.; 108; 1963; 275--292.

\bg\bg\hfill Math\'ematiques, laboratoire SDAD, FRE 2271
CNRS

\hfill Universit\'e Campus II, BP 5186, 14~032 Caen, France

\hfill dehornoy@math.unicaen.fr

\hfill http://www.math.unicaen.fr/$\sim$dehornoy/

\bye